\documentclass[a4paper,11pt,oneside]{article}

\usepackage{amsmath,amsthm,amsfonts,amssymb,amscd}
\usepackage{mathtools}
\usepackage{float}
\usepackage{subfig} 
\usepackage{enumerate}
\usepackage{booktabs}
\usepackage[noend]{algpseudocode}
\usepackage[normalem]{ulem}
\usepackage{graphicx,bm,xcolor}
\usepackage{algorithm}
\usepackage{algpseudocode}
\usepackage{hyperref} 
\usepackage{caption}
\usepackage[T1]{fontenc}
\usepackage{t1enc}
\usepackage[polish,german,english]{babel}
\usepackage{verbatim}
\usepackage{diagbox}
\usepackage{colortbl}
\usepackage{siunitx}
\usepackage{mathdots} 
\usepackage{bbm}

\usepackage{pgf,tikz,pgfplots}
\pgfplotsset{compat=1.18}
\usepackage{mathrsfs}
\usetikzlibrary{arrows.meta}
\usetikzlibrary{shapes.misc}
\usetikzlibrary{positioning}
\usetikzlibrary{decorations.pathreplacing} 

\usepackage{url}
\usepackage{authblk}

\algnewcommand{\algorithmicgoto}{\textbf{go to}}
\algnewcommand{\Goto}[1]{\algorithmicgoto~\ref{#1}}

\usepackage[colorinlistoftodos,prependcaption,textsize=scriptsize,color=red!40!white]{todonotes}\usepackage{totcount}

\definecolor{orange}{RGB}{230, 159, 0}
\definecolor{skyblue}{RGB}{86, 180, 233}
\definecolor{yellow}{RGB}{240, 228, 66}
\definecolor{blue}{RGB}{0, 114, 178}
\definecolor{vermillion}{RGB}{213, 94, 0}

\usepackage{listings}
\definecolor{codegreen}{rgb}{0,0.6,0}
\definecolor{codegray}{rgb}{0.5,0.5,0.5}
\definecolor{codepurple}{rgb}{0.58,0,0.82}
\definecolor{backcolour}{rgb}{0.95,0.95,0.96}
\definecolor{lightgray}{rgb}{1,1,1}

\lstdefinestyle{mystyle}{ 
    commentstyle=\color{codegreen},
    keywordstyle=\color{blue},
    numberstyle=\tiny\color{codegray},
    stringstyle=\color{codepurple},
    basicstyle=\ttfamily\footnotesize,
    breakatwhitespace=false,         
    breaklines=true,                 
    captionpos=b,                    
    keepspaces=true,                 
    numbers=left,                    
    numbersep=5pt,                  
    showspaces=false,                
    showstringspaces=false,
    showtabs=false,                  
    tabsize=2,
    xleftmargin=3pt
}
 
\usepackage[top=2cm, bottom=2cm, left=2cm, right=2cm]{geometry}

\theoremstyle{plain} 

\theoremstyle{definition} %

\theoremstyle{remark} %

\DeclareMathOperator*{\true}{\text{true}}
\DeclareMathOperator*{\guess}{\text{guess}}
\DeclareMathOperator*{\NN}{\text{NN}}

\definecolor{brewer1}{HTML}{A6CEE3}
\definecolor{brewer2}{HTML}{1F78B4}
\definecolor{brewer3}{HTML}{B2DF8A}
\definecolor{brewer4}{HTML}{33A02C}

\begin{document}

\title{
Optimal control of partial differential equations in PyTorch using automatic differentiation and neural network surrogates}
\author[1]{Denis Khimin}
\author[1]{Julian Roth}
\author[2]{Alexander Henkes}
\author[1,3]{Thomas Wick}

\affil[1]{Leibniz Universit\"at Hannover, Institut f\"ur Angewandte
  Mathematik, \newline AG Wissenschaftliches Rechnen, Welfengarten 1, 30167 Hannover, Germany}
\affil[2]{Computational Mechanics Group, ETH Zürich, Switzerland}
\affil[3]{Universit\'e Paris-Saclay, CentraleSupélec, ENS Paris-Saclay, CNRS, LMPS - Laboratoire de Mécanique Paris-Saclay,
91190 Gif-sur-Yvette, France}

\date{}
\maketitle

\setcounter{page}{1}


\begin{abstract}
The take-home message of this paper is that solving optimal control problems can be computationally straightforward, provided that differentiable partial differential equation (PDE) solvers are available.
Although this might seem to be a strong limitation and the development of differentiable PDE solvers might seem arduous, for many problems this is not the case.
In particular, for linear partial differential equations, they are equivalent to a linear equation system.
Therefore it is just sufficient to be able to solve linear equation systems in an automatic differentiation-capable library like PyTorch and be able to differentiate through the linear solver.
Using open-source libraries, like \texttt{torch\_sparse\_solve}, we have easy to use differentiable direct solvers (here: KLU) at hand, which makes solving linear PDEs straightforward.
For nonlinear PDEs, the approach above might be not sufficient, since we do not have an equivalent linear equation system for the underlying PDE. In this case, we can then use the library \texttt{torch-fenics}, which enables us to define PDEs using the finite element library \texttt{FEniCS} and then differentiate through this problem in PyTorch. We complement the proposed methodology with various optimal control problems constrained by stationary and instationary PDEs, linear and nonlinear PDEs (from Poisson to fluid-structure interaction), space-time and time-stepping formulations, parameter estimation, right-hand side control, initial condition control and boundary condition control, finite difference and finite element discretizations and neural network surrogates.  
All developments of this work are accompanied with the respective source codes published on GitHub.

\textbf{Highlights}
\begin{itemize}
    \item automatic differentiation for linear systems
    \item automatic differentiation for PDE-constrained optimal control
    \item parameter estimation
    \item optimal control (right hand side/boundary/initial value control)
    \item finite difference and finite element discretization
    \item neural network-based optimal control
\end{itemize}
\end{abstract}

\section{Introduction}
\label{sec:intro}

In many real-world problems and practical applications that are governed by PDEs, we only have limited knowledge of the PDEs and its initial condition, rather than full information.
Nevertheless, we have some measurements of the true solution and aim to infer the material parameters (parameter estimation), the boundary conditions (boundary control), the initial value (initial condition control), or, more generally, to minimize a specific quantity of interest, e.g. drag or lift coefficient.
These problems belong to the field of optimal control with partial differential equations \cite{lions1971, hinze2008optimization, Troltzsch2010}.
Often, the barrier to entry in this field is high for people with little background in numerical optimization. The PDE-constrained problem is transformed into an unconstrained optimization problem by the method of Lagrange multipliers \cite{LagrangeMultiplier}.
To solve this unconstrained problem, one has to write down its first order optimality conditions which are given by the Karush-Kuhn-Tucker (KKT) conditions \cite{karush1939, Kuhn1951}. This requires the derivation and implementation of an auxiliary adjoint equation (e.g.,~\cite{BeMeVe07}), and in particular, for higher order optimization schemes and highly nonlinear PDEs, the derivation, theory, together with the implementation and verification of these auxiliary problems can be cumbersome as in fluid-structure interaction~\cite{RiWi2013fsi,FaiMeiVex16,Feppon2019,Failer2021,Haubner_2020,FEPPON2020109574,WiWo21,haubner2024numericalmethodsshapeoptimal}.

As our focus in this work is on automatic differentiation, and consequently on efficient software implementations, we list some numerical optimization libraries in the following.
Several software libraries have been developed over the years. Here is a brief overview of some of them: RoDoBo \cite{becker2005rodobo} is a C++ library designed for optimization problems with PDE constraints. 
DOpElib \cite{DOpElib} is an open source library (C++) focused on PDEs and optimal control problems, particularly utilizing the reduced approach.
Dolfin, part of the FEniCS project \cite{fenicsproject}, provides an environment for solving PDEs with a Python and C++ interface. 
SNOPT (Sparse Nonlinear OPTimizer) \cite{Gill2005} solves large-scale nonlinear optimization problems (NLPs) utilizing the sequential quadratic programming method and excelling in sparse structures. Another robust and efficient library is IPOPT 
\cite{Wachter2005}, which primarily employs interior-point methods and is widely used in various domains like logistics, energy, and finance.

Automatic differentiation \cite{GriewankWalther2008, Naumann2011} is an easy-to-use solution to avoid the manual calculation of derivatives and auxiliary problems (e.g., adjoint problem, tangent problem, adjoint Hessian). As the name suggests, the derivatives are computed automatically which reduces the programmer's workload. Behind the scenes, the entire code, e.g. the PDE solver, is represented as a computational graph consisting of differentiable blocks allowing for the computation of the derivatives by applying the chain rule.
The gradients can be computed with the chain rule by traversing the computational graph either forward (forward mode automatic differentiation) or backward (reverse mode automatic differentiation). In this work, we employ the latter. Consequently, our optimal control solution consists of two different phases: the forward pass and the backward pass. As a user, we only want to implement the forward pass, which boils down to solving the original problem. The backward pass, which is required for the computation of the derivatives, should be done automatically. This, the entire automatic-differentation driven optimal control problem requires the following steps: implement the PDE solver, implement your quantity of interest or loss function depending on the observation data, wrap the optimization loop around the PDE solution and loss computation.

Regarding the differentiable PDE solver required for the forward pass, we rely on differentiable direct solvers \cite{LaporteBlog, LaporteCode} and finite element libraries with automatic differentiation capabilities \cite{Mitusch2019, BarkmanCode}.
More concretely, since linear PDEs are equivalent to solving a linear equation system, the forward pass only requires the solution of that linear equation system. This can be dune using a differentiable direct solver based on the KLU method \cite{KLU2010}, developed by Floris Laporte in the Python package \texttt{torch\_sparse\_solve} \cite{LaporteCode, LaporteBlog}.
For nonlinear PDEs having a differentiable linear solver is not enough. Instead we use automated differentiation combined with finite element method libraries, e.g. in \texttt{dolfin-adjoint} \cite{Mitusch2019, farrell2013, funke2013} for the \texttt{FEniCS} library \cite{fenicsproject}.
The backward pass, i.e. the adjoint equation, is then derived automatically by the finite element library, e.g by \texttt{torch-fenics} \cite{BarkmanCode} which is Patrik Barkman's wrapper around \texttt{dolfin-adjoint} that converts \texttt{FEniCS} data structures to \texttt{PyTorch} data structures. 
Due to the PyTorch wrappers, it is then also possible to use neural networks in the optimal control workflow as shown in Section \ref{sec:nn_optimal_control}.

The outline of this paper is as follows. 
In Section \ref{sec:automatic_differentiation}, we give a short introduction to automatic differentiation, with applications to linear PDEs in Section \ref{sec:automatic_differentiation_linear_system} and general nonlinear PDEs in Section \ref{sec:automatic_differentiation_PDE_solver}.
Next, in Section \ref{sec:numerical_tests}, we present nine different optimal control problems with finite difference and finite element discretizations with examples ranging from the Poisson problem to nonlinear fluid-structure interaction, including also an example for neural network controls.
Our work is summarized in Section \ref{sec:conclusion}.

\section{Automatic differentiation}
\label{sec:automatic_differentiation}

Automatic differentiation \cite{GriewankWalther2008, Naumann2011} is a set of techniques used to efficiently and accurately calculate derivatives of functions that are implemented as computer programs, e.g. codes for the solution of partial differential equations.
The main idea of automatic differentiation is the interpretation of the code as a computational graph consisting of differentiable blocks such that the derivatives can be applied through the successive application of the chain rule. More concretely, in this work we employ reverse mode automatic differentiation, which is also used for the backpropagation of gradients in neural networks \cite{LeCun2012} and is being used in many finite element libraries, e.g. in \texttt{dolfin-adjoint} \cite{Mitusch2019, farrell2013, funke2013} for the \texttt{FEniCS} library \cite{fenicsproject}. Reverse mode automatic differentiation entails two phases: the forward pass and the backward pass.
In the forward pass, the original problem is being solved; for example, the linear system is solved, the PDE is solved (forward in time), or the neural network makes its prediction.
In the backward pass, starting at the output from the forward pass, gradients (adjoints) are being passed backwards through the computational graph, e.g. an auxiliary linear system or an auxiliary PDE is solved (backward in time), or gradients are computed for each layer starting at the output layer and ending at the input layer. 
In the next step, we will describe in greater detail how the reverse mode automatic differentiation works for linear systems and PDE solvers based on the finite element method.

\subsection{Automatic differentiation for linear systems}
\label{sec:automatic_differentiation_linear_system}
\begin{figure}[H]
    \begin{center}
    \begin{tikzpicture}[>=stealth, node distance=2cm, on grid, auto]

      \draw[->, thick, blue] (0,1.15) -- node[above] {forward} (9,1.15);
      \draw[->, thick, purple] (9,1) -- node[below] {backward} (0,1);
      
      \node[circle, draw=black, thick, minimum size=1cm] (A) at (0,0) {A};
      \node[circle, draw=black, thick, minimum size=1cm, below of=A, node distance=2cm] (B) {b};
      \node[rectangle, draw=black, thick, minimum height=1cm, minimum width=2.5cm, right of=A, node distance=3cm, yshift=-1cm] (solve) {solve(A, b)};
      \node[circle, draw=black, thick, minimum size=1cm, right of=solve, node distance=3cm] (X) {x};
      \node[circle, draw=black, thick, minimum size=1cm, right of=X, node distance=3cm] (J) {J};
    
      \draw[->, thick] (A) -- (solve);
      \draw[->, thick] (B) -- (solve);
      \draw[->, thick] (solve) -- (X);
      \draw[->, thick] (X) -- (J);
    
    \end{tikzpicture}
    \caption{Computation graph for the solution of a linear equation system.}
    \label{fig:linear_system}
    \end{center}
\end{figure}
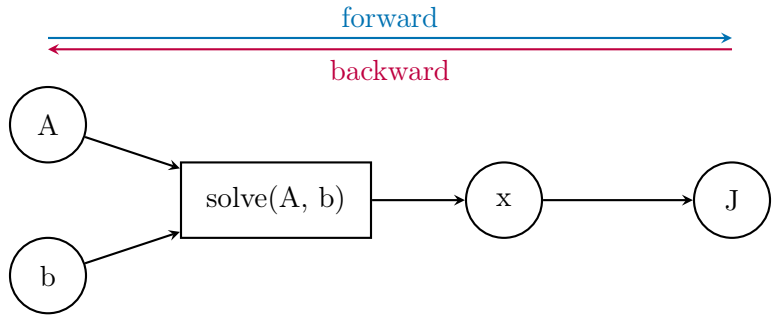

We consider the linear system
\begin{align*}
    Ax = b
\end{align*}
with $A \in \mathbb{R}^{n \times n}, b \in \mathbb{R}^n, n \in \mathbb{N}$.
The forward pass takes $A$ and $b$ as input, computes the solution of the linear system
\begin{align*}
    x = A^{-1}b = \operatorname{solve}(A, b)
\end{align*}
and then computes some scalar-valued quantity of interest (loss function) $J = J(x)$, e.g. $J(x) = \|x - \tilde{x}\|_2$ with some reference solution $\tilde{x} \in \mathbb{R}^n$. The corresponding computation graph is shown in Figure \ref{fig:linear_system}.

For the backward pass, we know the gradient of the quantity of interest $J$ w.r.t. the solution of the linear system $x$, i.e. we know $\frac{\partial J}{\partial x}$. From this, we can compute $\frac{\partial J}{\partial A}$ and $\frac{\partial J}{\partial b}$ by the relations
\begin{align*}
    \frac{\partial J}{\partial b} &= A^{-T} \frac{\partial J}{\partial x} = \operatorname{solve}\left(A^T, \frac{\partial J}{\partial x}\right), \\
    \frac{\partial J}{\partial A} &= - \frac{\partial J}{\partial b} \otimes x.
\end{align*}
Now, that we know how to differentiate through a single linear system, we can use this as one of many differentiable blocks in a larger computational graph, e.g. we can differentiate through time stepping schemes where in each time step a linear system must be solved. For the implementation, a differentiable direct solver based on the KLU method \cite{KLU2010} has been developed by Floris Laporte in the Python package \texttt{torch\_sparse\_solve} \cite{LaporteCode, LaporteBlog}. Instead of differentiating through a direct solver, one could also differentiate through an iterative solver of the linear system \cite{hovland2024differentiatinglinearsolvers}. 

The abstract workflow for solving optimal control problems based on linear systems using automatic differentiation is given in Listing \ref{lst:pseudo_linear_system}. To solve the optimization problem we need to follow the steps below.
First, we need to initialize the parameter guess \lstinline{param_guess} and the true parameter \lstinline{param_true}. 
Next, we determine the solutions of the linear systems depending on these parameters, i.e. the solution \lstinline{x_guess} where the matrix $A$ and the right hand side $b$ may depend on the parameter guess \lstinline{param_guess} and the solution \lstinline{x_true} where the matrix $A$ and the right hand side $b$ may depend on the true parameter \lstinline{param_true}.
Then, after choosing a suitable optimizer, e.g. Rprop \cite{Rprop1993}, and an appropriate learning rate, e.g. 0.1, we can continue with the optimization loop.
Note that the choice of the hyperparameters, like optimizer type, learning rate, number of optimization steps, is crucial for the success of the optimization and strongly depends on the problem at hand. The optimization loop itself should look very familiar to the optimization loops that we are accustomed to from neural network training. Before the computation, the gradients need to be zeroed out. Thereafter, we compute our prediction by solving the linear system given the current parameter guess.
This prediction is then plugged into the loss function $J$ that evaluates how far our solution guess \lstinline{x_guess} deviates from the true solution \lstinline{x_true}.
Having completed the forward pass, we can now backpropagate the gradients and perform an optimization step.
We repeat the whole process of alternating forward and backward runs, until the target tolerance in the goal functional or the number of maximum iterations is reached.

\lstset{language=Python, caption=Pseudocode for automatic differentiation for linear systems, label={lst:pseudo_linear_system}}
\begin{lstlisting}[mathescape=true]
import torch
from torch_sparse_solve import solve

param_true, param_guess = ..., ...
# prepare the true solution
x_true = solve(A(param_true), b(param_true))
# prepare the initial guess
x_guess = solve(A(param_guess), b(param_guess))
# prepare optimization
iter, MAX_ITER = 0, 100
optimizer = torch.optim.Rprop([param_guess], lr=0.1)

# optimize the parameters as long as x_true and x_guess are not close enough
while J(x_true, x_guess) > 1e-6 and iter < MAX_ITER:
    iter += 1
    print(f"Iteration {iter}: Loss = {J(x_true, x_guess)}")
    # zero the gradients
    optimizer.zero_grad()
    # solve the linear system
    x_guess = solve(A(param_guess), b(param_guess))
    # compute the loss
    loss = J(x_true, x_guess)
    # backpropagate
    loss.backward()
    # update the parameter
    optimizer.step()
\end{lstlisting}

\subsection{Automatic differentiation for FEM-based PDE solvers}
\label{sec:automatic_differentiation_PDE_solver}
In some cases having a differentiable direct solver is not enough. In contrast to linear problems, nonlinear problems cannot be written as a linear system (without linearization).
Therefore, we present a more general approach. 
Let $Q$ be the control space and $V$ the function space for the PDE solution. Then, an abstract PDE-constrained optimization can be formulated as
\begin{align*}
    \min_{q, u} J(q, u) \qquad \text{s.t.} \qquad A(q, u)(\psi) = 0 \quad \forall \psi \in V.
\end{align*}
The goal is to minimize the quantity of interest $J$ that depends on the control variable $q\in Q$ and the solution $u\in V$ of the PDE 
\[
A(q, u)(\psi) = 0 \quad \forall \psi \in V.
\]
Here, $A(q, u)(\psi)$ represents the weak form of the PDE. All problem statements from Section~\ref{sec:numerical_tests} can be framed into such a weak form by employing appropriate function spaces $V$, where possibly 
$V$ need to be further refined for the trial and test functions depending 
on the boundary conditions.

The Lagrange functional with the Lagrange multiplier,
i.e., the adjoint $z\in V$ (in general $z\in V^*$) for this optimization problem is given by
\begin{align*}
    \mathcal{L}(q, u, z) := J(q, u) - A(q, u)(z).
\end{align*}
The first order optimality conditions are given by the Karush-Kuhn-Tucker (KKT) conditions \cite{karush1939, Kuhn1951}
\begin{align}
    \mathcal{L}_u^{\prime}(q, u, z)(\delta u) = J_u^{\prime}(q, u)(\delta u) - A_u^{\prime}(q, u)(\delta u, z) = 0 \quad \forall \delta u \in V, \tag{adjoint equation} \\
    \mathcal{L}_q^{\prime}(q, u, z)(\delta q) = J_q^{\prime}(q, u)(\delta q) - A_q^{\prime}(q, u)(\delta q, z) = 0 \quad \forall \delta q \in Q, \tag{gradient equation} \\
    \mathcal{L}_z^{\prime}(q, u, z)(\delta z) = - A(q, u)(\delta z) = 0 \quad \forall \delta z \in V. \tag{state equation}
\end{align}
Herein we assume that the partial derivatives of 
$J$ and $A$ exist in a weak sense. 
First of all, we need to solve the state equation to obtain a solution $u \in V$ of the original PDE.
Then, we have to solve the adjoint equation to obtain a solution $z \in V$ of the auxiliary problem (adjoint) that measures the sensitivity of the PDE w.r.t. the goal functional $J$.
Plugging the primal and dual solutions into the gradient equation, we can then compute the gradient of the goal functional $J$ w.r.t. the control $q \in Q$.
Luckily we do not need to manually compute all the derivatives $J_{\Box}^{\prime}$ and $A_{\Box}^{\prime}$ for $\Box \in \{u, q\}$, but can use automated differentiation within finite element method libraries, e.g. in \texttt{dolfin-adjoint} \cite{Mitusch2019, farrell2013, funke2013} for the \texttt{FEniCS} library \cite{fenicsproject}.
By this, it is sufficient to define the variational formulation of the PDE $A(q, u)(\psi) = 0$ and the goal functional $J(q, u)$.
For the numerical examples, we use \texttt{torch-fenics} \cite{BarkmanCode}, which is Patrik Barkman's wrapper around \texttt{dolfin-adjoint} that converts \texttt{FEniCS} data structures to \texttt{PyTorch} data structures. Consequently, it is also easier to incorporate neural network surrogates into the optimal control workflow as shown in Section \ref{sec:nn_optimal_control}.

\subsection{Example of implementation details}
We will discuss the implementational details using the thermal fin example from Section \ref{sec:2d_thermal_fin}. To use \texttt{torch-fenics} \cite{BarkmanCode}, we need to implement two things: the FEniCS-based PDE solver class (Listing \ref{lst:fem_pde_solver_fenics_class}) and the main function that contains the optimization loop (Listing \ref{lst:fem_pde_solver_main}).
First, let us take a closer look at the PDE solver in Listing \ref{lst:fem_pde_solver_fenics_class}. After importing some libraries, we define a new class, here called \lstinline{Poisson}, that implements the abstract base class \lstinline{torch_fenics.FEniCSModule}. This abstract base class contains wrappers that handle the conversion between PyTorch, NumPy and FEniCS datatypes, as well as the solution of the backward problem.
The user then only needs to write a class function \lstinline{solve} that takes the parameters as input and solves the PDE, and a class function \lstinline{input_templates} that defines the input shape of the parameters that are being passed to \lstinline{solve}.
The \lstinline{solve} method follows standard FEniCS \cite{fenicsproject} syntax and solves a Poisson problem with different heat conductivies for the respective subdomains as well as Neumann and Robin boundary conditions (cf. Section \ref{sec:2d_thermal_fin}).

\lstset{language=Python, caption=Truncated FEniCS PDE class from Section \ref{sec:2d_thermal_fin} for automatic differentiation for FEM-based PDE solvers, label={lst:fem_pde_solver_fenics_class}}
\begin{lstlisting}[mathescape=true]
import torch
from fenics import *
from fenics_adjoint import *
import torch_fenics

class Poisson(torch_fenics.FEniCSModule):
    def __init__(self):
        ...
        
    def solve(self, mu):
        # Create trial and test functions
        u = TrialFunction(self.V)
        v = TestFunction(self.V)

        # Construct bilinear form:
        # * subdomain integrals with different heat conductivities mu[i] = k_i
        a = mu[0] * inner(grad(u), grad(v)) * self.dx(1) + mu[1] * inner(grad(u), grad(v)) * self.dx(2) + mu[2] * inner(grad(u), grad(v)) * self.dx(3) + mu[3] * inner(grad(u), grad(v)) * self.dx(4) + mu[4] * inner(grad(u), grad(v)) * self.dx(5)
        # * boundary integral for Robin boundary condition with heat transfer coefficient mu[5] = Bi
        a += mu[5] * u * v * self.ds(2)
        # Construct linear form
        L = Constant(1.) * v * self.ds(1)

        # Solve the Poisson equation
        u = Function(self.V)
        solve(a == L, u)

        # Return the solution
        return u

    def input_templates(self):
        # Declare templates for the inputs to Poisson.solve
        return Constant((0, 0, 0, 0, 0, 0))
\end{lstlisting}

The optimization loop for \texttt{torch-fenics} \cite{BarkmanCode} is shown in Listing \ref{lst:fem_pde_solver_main}. 
We observe that this code is very similar to the optimization loop in Listing \ref{lst:pseudo_linear_system}, which proves the point that optimal control is easy to implement as we just need to implement the forward solver and tune some optimization parameters, whereas the rest of the code remains unchanged. 
We want to point out two novelties compared to Listing \ref{lst:pseudo_linear_system} which are specific to the optimal control problem from Section \ref{sec:2d_thermal_fin}.
It is common to use regularization in the loss function, which is being included here through the term $+ 0.1 \left\|\frac{\mu^{\guess}-\mu^{\operatorname{ref}}}{\mu^{\operatorname{ref}}}\right\|_2$.
Moreover, we need to ensure that the parameters $\mu$ are bounded by $[0.1, 10]^5 \times [0.01, 1]$ which is being done by \lstinline{mu_guess.clamp_} after the optimization step. This projects the parameters $\mu$ on the boundary of the parameter domain, if the constraints are not satisfied.

\lstset{language=Python, caption=Truncated main function from Section \ref{sec:2d_thermal_fin} for automatic differentiation for FEM-based PDE solvers, label={lst:fem_pde_solver_main}}
\begin{lstlisting}[mathescape=true]
# Construct the FEniCS model
poisson = Poisson()
# mu = [k0, k1, k2, k3, k4, Bi] -> parameters in PDE which are to be learned
mu_true = torch.tensor([[0.1, 8.37317446377103, 6.572276066240383, 0.46651735398635275, 1.8835410659596712, 0.01]], dtype=torch.float64)
# get the true solution of the Poisson equation
u_true = poisson(mu_true)
# get a reference mu for regularization
mu_reference = torch.tensor([[1., 1., 1., 1., 1., 0.1]], dtype=torch.float64)
# perform optimization of mu_guess
mu_guess = torch.tensor(0.5 * torch.ones(1, 6, dtype=torch.float64), requires_grad=True) # = [0.5, 0.5, 0.5, 0.5, 0.5, 0.5]
u_guess = poisson(mu_guess)

# prepare optimization
iter, MAX_ITER = 0, 100
optimizer = torch.optim.Rprop([mu_guess], lr=0.01)
print("Optimizing the parameters in the thermal fin problem...")
print(f"Number of parameters: {mu_guess.numel()}")

# optimize the parameters as long as u_true and u_guess are not close enough
while torch.norm(u_true - u_guess) > 1e-6 and iter < MAX_ITER:
    iter += 1
    print(f"Iteration {iter}: Loss = {torch.norm(u_true - u_guess)}")
    # zero the gradients
    optimizer.zero_grad()
    # solve the Poisson equation
    u_guess = poisson(mu_guess)
    # compute the loss
    loss = torch.norm(u_true - u_guess) + 0.1 * torch.norm((mu_guess - mu_reference) / mu_reference)
    # backpropagate
    loss.backward()
    # update the parameters
    optimizer.step()
    # apply constraints
    with torch.no_grad():
        mu_guess.clamp_(
            min=torch.tensor([0.1, 0.1, 0.1, 0.1, 0.1, 0.01], dtype=torch.float64), 
            max=torch.tensor([10., 10., 10., 10., 10., 1.], dtype=torch.float64)
        )
\end{lstlisting}

\subsection{Automatic differentiation for PDE solvers using neural network surrogates}
\label{sec:automatic_differentiation_using_NN}
To combine neural networks with optimal control, we focus on feedforward neural networks $\operatorname{NN}: \bar{\Omega} \subset \mathbb{R}^2 \rightarrow \mathbb{R}$ in this work.
The neural networks can be expressed as
\begin{align*}
    \operatorname{NN}(x) = T^{(L)} \circ \sigma \circ T^{(L-1)} \circ \cdots \circ \sigma \circ T^{(1)}(x),
\end{align*}
where $T^{(i)}: \mathbb{R}^{n_{i-1}} \rightarrow \mathbb{R}^{n_{i}},y \mapsto  W^{(i)} y\, +\, b^{(i)}$ are affine transformations for $1 \leq i \leq L$, with weight matrices $W^{(i)} \in \mathbb{R}^{n_i \times n_{i-1}}$ and bias vectors $b^{(i)} \in \mathbb{R}^{n_i}$. Here $n_{i}$ denotes the number of neurons in the $i$--th layer with $n_0 = 2$ and $n_L = 1$, $\sigma: \mathbb{R} \rightarrow \mathbb{R}$ is a nonlinear activation function, which is the sigmoid function throughout this work.
The main idea is to replace the control vector by a neural network which is being evaluated at the degrees of freedom of the control variable $q$. The code for the optimization loop from Listing \ref{lst:fem_pde_solver_main} mostly remains the same. However, we have to feed the prediction of the neural network as a parameter input to the PDE solver. An excerpt from the code for Section \ref{sec:nn_optimal_control} is shown in Listing \ref{lst:pde_solver_nn}.

\lstset{language=Python, caption=Truncated optimization loop from Section \ref{sec:nn_optimal_control} for automatic differentiation for PDE solvers using neural network surrogates, label={lst:pde_solver_nn}}
\begin{lstlisting}[mathescape=true]
# Zero the gradients
optimizer.zero_grad()
# Forward pass
kappa = nn(input).T
# solve the Poisson equation
u_guess = poisson(kappa)
# compute the loss
loss = torch.norm(u_true - u_guess) 
# Backward pass
loss.backward()
# Update the parameters
optimizer.step()
\end{lstlisting}

We observe that the neural network based code in Listing \ref{lst:pde_solver_nn} is almost identical with the code in Listing \ref{lst:fem_pde_solver_main} and only the line \lstinline{kappa = nn(input).T} is new. Therein, the neural network \lstinline{nn} is being evaluated at the coordinates of the degrees of freedom \lstinline{input}, which in turn yields the prediction for the parameter $\kappa$ that is being plugged into the PDE.

\section{Numerical tests}
\label{sec:numerical_tests}

We perform nine computational experiments for various problems. 
Therein, the first four examples are rather from an educational viewpoint, while the following five numerical tests address more challenging PDEs such as incompressible flow, fluid-structure interaction, where multiple parameters must be estimated (i.e., learned).
In our final example we aim to find a 
neural network surrogate for a Poisson diffusion coefficient.

For the first four numerical tests, we perform computations using the finite difference method for the Poisson problem and the heat equation. 
We use a one-dimensional Poisson problem with a scalar-valued and a vector-valued force function that we want to recover from the true solution.

For the 1+1D heat equation, we use a space-time formulation and a time-stepping formulation to infer the initial condition (and right hand side).
For the next five numerical tests, we perform computations using the finite element method.
We estimate the material parameters for a 2D Poisson problem, infer the initial condition of a 2+1D nonlinear heat equation, minimize the drag through boundary control of the Navier-Stokes conditions and estimate the parameters of a stationary fluid-structure interaction problem. Finally, we demonstrate how neural networks can be combined with optimal control of a 2D Poisson problem.

The first four numerical tests have been implemented with \texttt{torch\_sparse\_solve} \cite{LaporteCode, LaporteBlog} and the remaining five numerical tests have been implemented with \texttt{torch-fenics} \cite{BarkmanCode}.

\subsection{1D Poisson with scalar-valued force}
\label{sec:1d_poisson_scalar}
\begin{figure}[H]
    \centering
    \includegraphics[width=12cm]{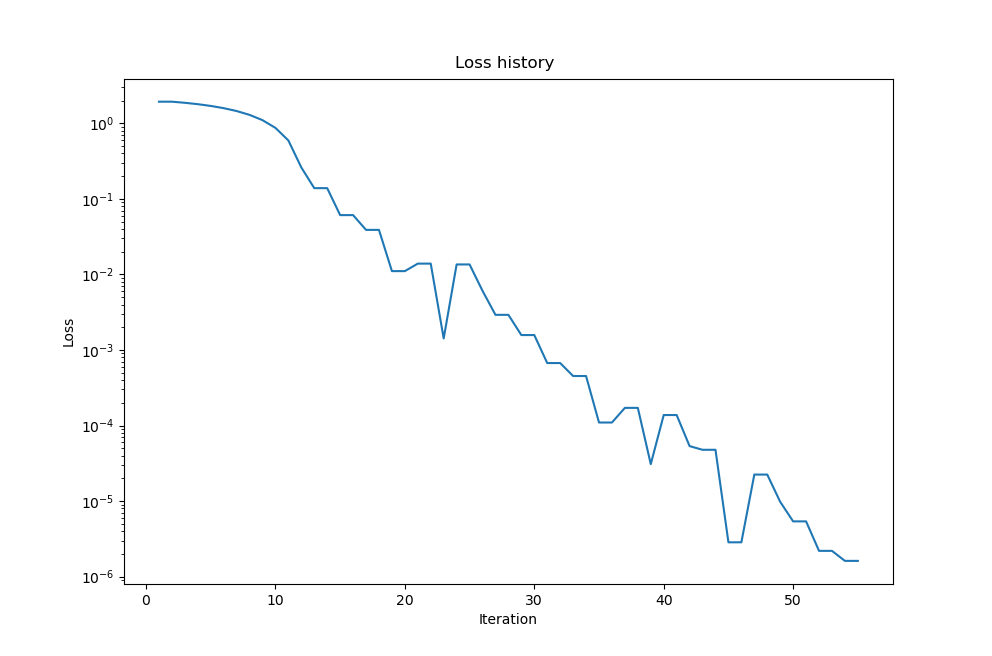}
    \caption{Example 1: Loss history of 1D Poisson with scalar-valued force.}\label{fig:loss_3_1}
\end{figure}

As a first numerical test, we consider the one-dimensional Poisson equation with a spatially-constant scalar-valued force. 
The goal of this example is to determine the gravitational force $f^{\guess}$ from given observations $U_h$ by comparing to some reference data $U_h(f^{true})$.
The strong form then reads: Find $u: (0,1) \rightarrow \mathbb{R}$ such that
\begin{align*}
    -\partial_x^2 u(x) &= f \quad \forall x \in (0,1), \\
    u(0) &= u(1) = 0.
\end{align*}
Using a uniform mesh with $x_{i+1} = x_i + h$ for $0 \leq i < n_h$, $x_0 = 0$, $x_{n_h} = 1$ and the spatial mesh size $h = \frac{1}{n_h} > 0$, we can use the central difference quotient
\begin{align*}
    \partial_x^2 u(x) \approx \frac{1}{h^2} \left(u(x+h) - 2u(x) + u(x-h) \right).
\end{align*}
Using this finite difference stencil with $u_i := u(x_i)$ and enforcing the boundary conditions $u_0 = u_{n_h} = 0$, we arrive at the linear system
\begin{align}\label{eq:1d_poisson_system}
    \frac{1}{h^2}
    \begin{bmatrix}
    2 & -1 & & & \\
    -1 & \ddots & \ddots & & \\
     & \ddots & \ddots & \ddots & \\
     & & \ddots & \ddots & -1 \\
     & & & -1 & 2
    \end{bmatrix}
    \begin{bmatrix}
        u_1 \\
        u_2 \\
        \vdots \\
        u_{n_h-2} \\
        u_{n_h-1}
    \end{bmatrix}
    &= 
    \begin{bmatrix}
        f \\
        f \\
        \vdots \\
        f \\
        f
    \end{bmatrix}, \\
    \Leftrightarrow \hspace{5cm} \qquad  
    K_h U_h &= F_h. \notag
\end{align}
For the numerical test, we choose $f^{\true} = -1$ and we start with an initial guess $f^{\guess} = 2$.
Using $n_h = 50$ spatial elements, the loss function
\begin{align*}
    J(f^{\guess}) := \|U_h(f^{true}) - U_h(f^{\guess}) \|_2,
\end{align*}
the Rprop (see Section~\ref{sec:automatic_differentiation})
optimizer with a learning rate
of $0.1$, after $55$ optimization steps we get the approximate result
\begin{align*}
    f^{\guess} \approx -0.999999583,
\end{align*}
which is close to the true solution up to a tolerance of $10^{-6}$. The loss history is shown in Figure \ref{fig:loss_3_1}.

\subsection{1D Poisson with vector-valued force}
\label{sec:1d_poisson_vector}
\begin{figure}[H]
    \centering
    \subfloat[Full]{
        \includegraphics[width=7cm]{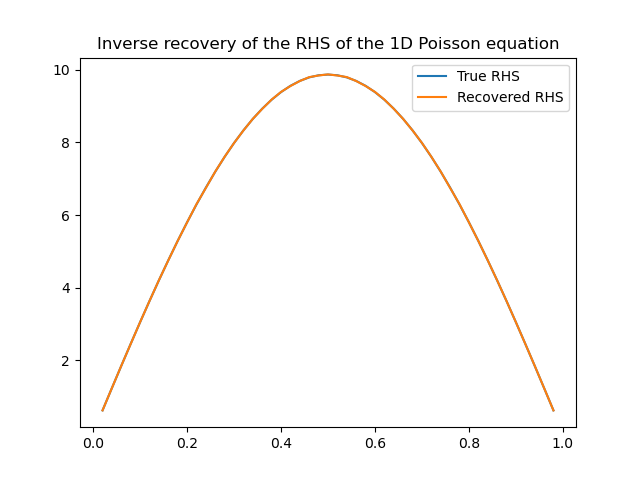}
    }%
    \quad
    \subfloat[Zoomed in]{
        \includegraphics[width=7cm]{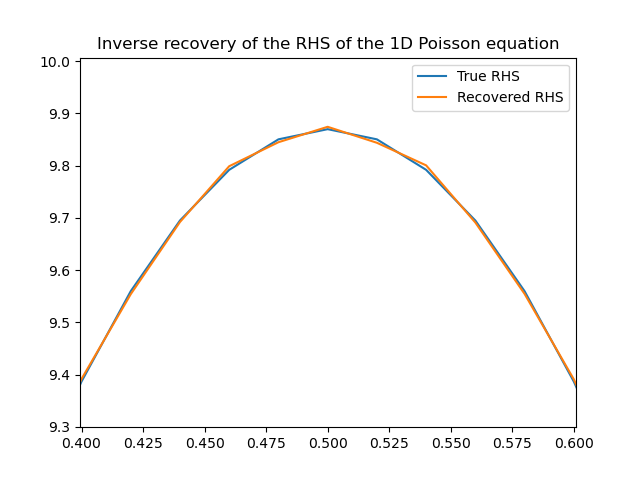}
    }
    \caption{Example 2: Recovered vector-valued force vector for a 1D Poisson problem.}\label{fig:recovered_rhs_1d_poisson}
\end{figure}

\begin{figure}[H]
    \centering
    \includegraphics[width=18cm]{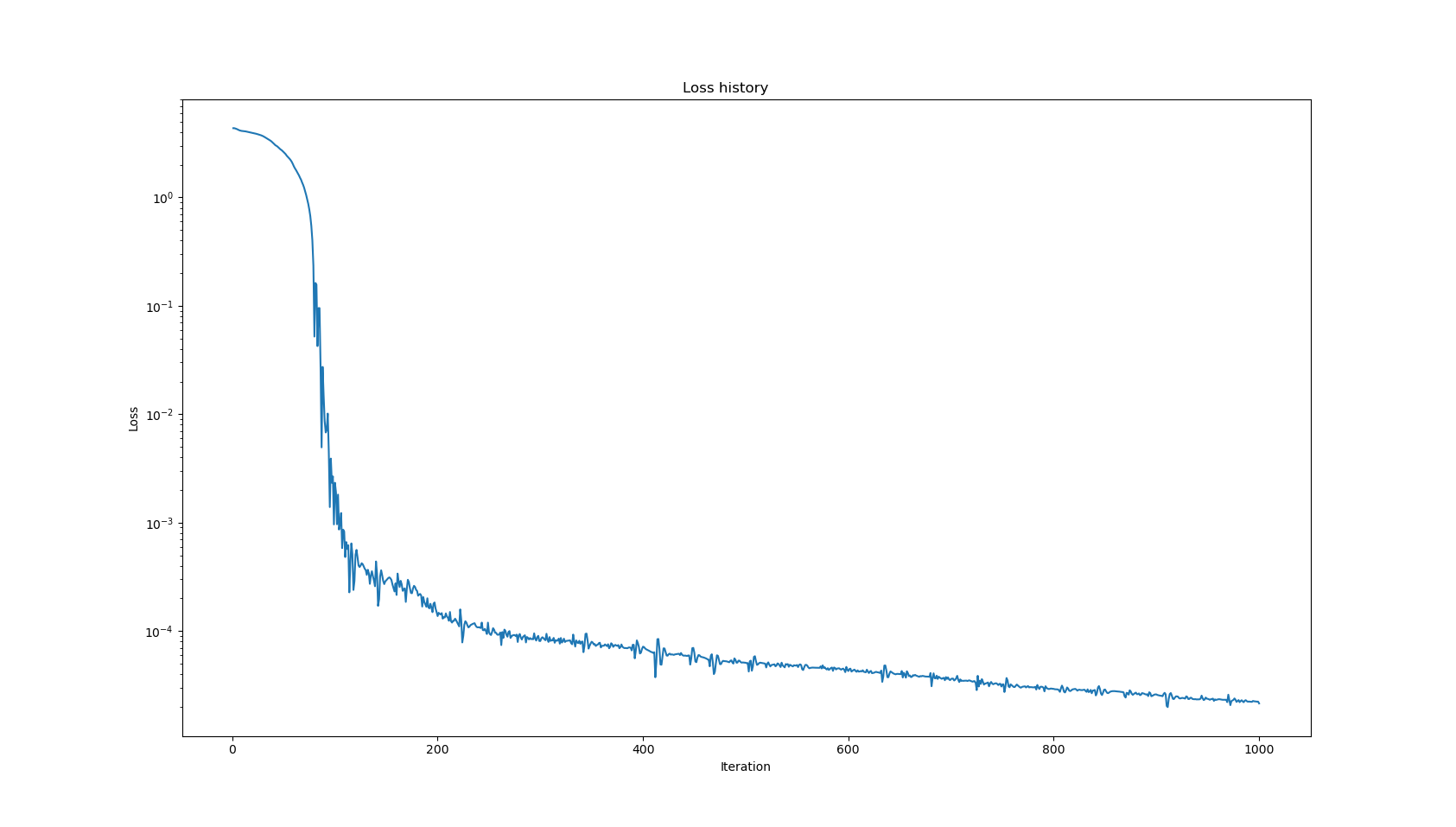}
    \caption{Example 2: Loss history of 1D Poisson with vector-valued force.}\label{fig:loss_3_2}
\end{figure}

In the second numerical test, we consider a very similar problem, where we replace the spatially constant force $f$ with a spatially-variable force $f(x)$.
In this example, we use
\begin{align*}
    f^{\true}(x) := \pi^2 \sin(\pi x).
\end{align*}
After a finite difference discretization,
we seek for a force vector
\begin{align*}
    F_h^{\guess} = \begin{pmatrix}f_1, f_2, \dots, f_{n_h-2}, f_{n_h-1}\end{pmatrix}^T \in \mathbb{R}^{n_h-1}
\end{align*}
with
\begin{align*}
    K_h U_h^{\true} \approx F_h^{\guess}.    
\end{align*}
In our loss function, we use Tikhonov regularization with $\alpha = 0.099$ 
and define it as
\begin{align*}
    J(F_h^{\guess}) := \|U_h(F_h^{true}) - U_h(F_h^{\guess}) \|_2 + \alpha \|F_h^{\guess} \|_2. 
\end{align*}
After $1000$ optimization steps, we reach a loss of $2.15 \cdot 10^{-5}$ and obtain the right hand side shown in Figure \ref{fig:recovered_rhs_1d_poisson}, where the loss history is shown in Figure \ref{fig:loss_3_2}.

\subsection{1+1D space-time heat equation with vector-valued force and initial condition}
\label{sec:1+1d_space_time_heat}
As the third numerical test, we consider a space-time discretization of the heat equation in one space-dimension, i.e. in 1+1D. 
The goal of this example is to determine $F^{\guess}$ from given space-time observations $U_{kh}$ by comparing it to some reference data $U_{kh}(f^{true})$.

The strong form reads: find $u: (0,1) \times (0, T) \rightarrow \mathbb{R}$ such that
\begin{align*}
    \partial_t u(x, t) - \partial_x^2 u(x,t) &= f(x,t), \qquad \forall (x,t) \in (0,1) \times (0, T), \\
    u(0, t) &= u(1, t) = 0, \quad \forall t \in (0, T), \\
    u(x, 0) &= u_0(x), \qquad \forall x \in (0, 1).
\end{align*}
Again, we use a uniform mesh in space with $x_{i+1} = x_i + h$ for $0 \leq i < n_h$, $x_0 = 0$, $x_{n_h} = 1$ and the spatial mesh size $h = \frac{1}{n_h} > 0$.
Similarly, we use a uniform mesh in time with $t_{j+1} = t_j + k$ for $0 \leq j < n_k$, $t_0 = 0$, $t_{n_k} = T$ and the temporal mesh size $k = \frac{T}{n_k} > 0$.
We discretize the derivatives with the finite differences stencils
\begin{align*}
    \partial_x^2 u(x, t) &\approx \frac{1}{h^2} \left(u(x+h, t) - 2u(x, t) + u(x-h, t) \right), \\
    \partial_t u(x, t) &\approx \frac{1}{k} \left(u(x, t) - u(x, t-k) \right).
\end{align*}
Using this finite difference stencil with $u_{i,j} := u(x_i, t_j)$ and enforcing the homogeneous boundary conditions $u_{0, j} = u_{n_h, j} = 0$ and the initial condition $u_{i, 0} = u_0(x_i)$, we arrive at the space-time linear system
\begin{align}\label{eq:heat_1+1D_space_time_system}
    \begin{bmatrix}
    \frac{1}{k}I_h + K_h & & & \\
    -\frac{1}{k}I_h & \ddots & & \\
    & \ddots & \ddots & \\
    & & \ddots & \frac{1}{k}I_h + K_h \\
    & & & -\frac{1}{k}I_h 
    \end{bmatrix}
    \begin{bmatrix}
        u_{\bullet, 1} \\
        u_{\bullet, 2} \\
        \vdots \\
        u_{\bullet, n_k-1} \\
        u_{\bullet, n_k}
    \end{bmatrix}
    &= 
    \begin{bmatrix}
        f_{\bullet, 1} + \frac{1}{k} u_{\bullet, 1}\\
        f_{\bullet, 2} \\
        \vdots \\
        f_{\bullet, n_k-1} \\
        f_{\bullet, n_k}
    \end{bmatrix}, \\
    \Leftrightarrow \hspace{1cm} \qquad  
    \underbrace{\left[I_k \otimes \left(\frac{1}{k}I_h + K_h\right) - \frac{1}{k}S_k \otimes I_h\right]}_{=: A_{kh}} U_{kh} &= F_{kh}. \notag
\end{align}
Here, we used the stiffness matrix $K_h$ from the Poisson problem in equation (\ref{eq:1d_poisson_system}), $I_k \in \mathbb{R}^{n_k \times n_k}$ and $I_h \in \mathbb{R}^{(n_h-1) \times (n_h-1)}$ are identity matrices and $S_k$ is the lower shift matrix with ones on the subdiagonal, i.e. $(S_k)_{i,j} = \delta_{i,j+1}$ for $1 \leq i,j \leq n_k$. By $A \otimes B$ we denote the Kronecker product of two matrices $A \in \mathbb{R}^{n_k \times n_k}$ and $B \in \mathbb{R}^{(n_h-1) \times (n_h-1)}$, which is given by
\begin{align*}
    A \otimes B = \begin{pmatrix}
        A_{1,1}B & \cdots & A_{1,n_k}B \\
        \vdots & \ddots & \vdots \\
        A_{n_k,1}B & \cdots & A_{n_k,n_k}B
    \end{pmatrix} \in \mathbb{R}^{n_k(n_h-1) \times n_k(n_h-1)}.
\end{align*}

\begin{figure}[H]
    \centering
    \subfloat[Recovered solution]{
        \includegraphics[width=7.5cm]{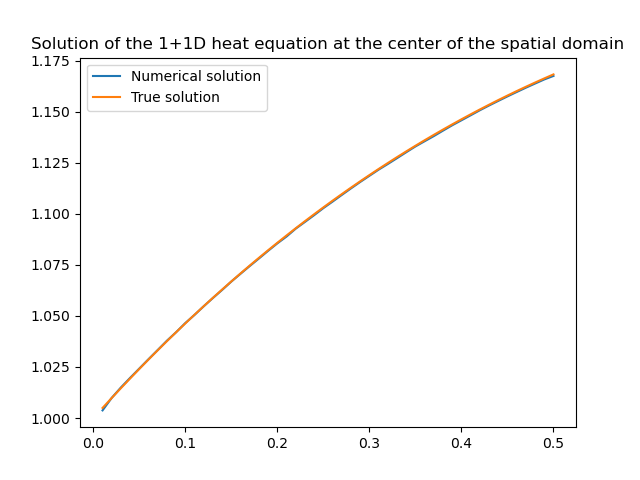}
    }%
    \subfloat[Loss history]{
        \includegraphics[width=9.75cm]{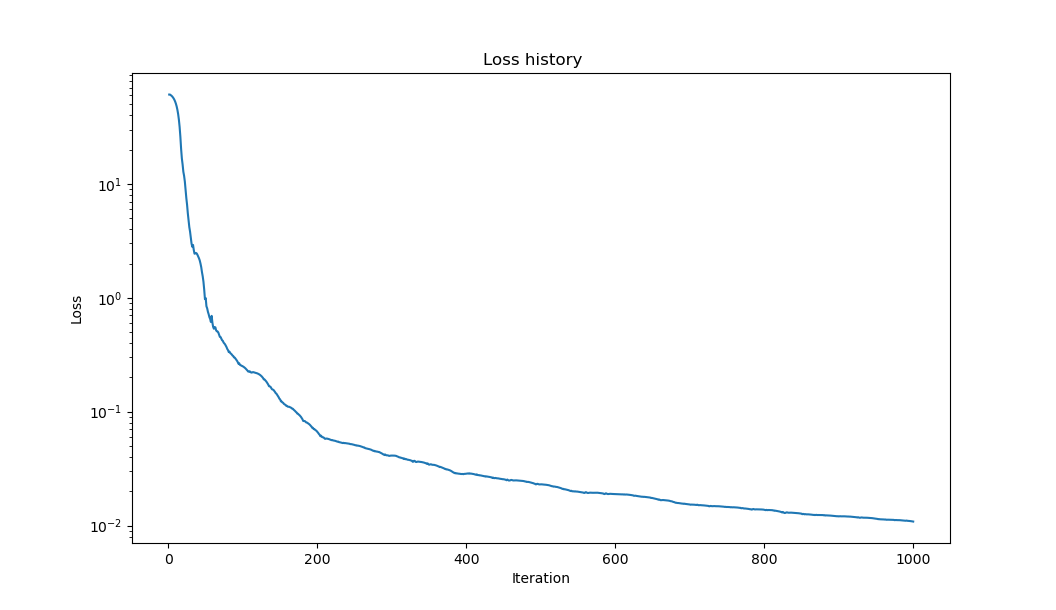}
    }%
    \caption{Example 3: Recovered solution of the space-time heat equation in 1+1D at $x = \frac{1}{2}$ (a) and
    loss history (b).}
    \label{fig:recovered_rhs_1+1d_heat_center}
    \label{fig:loss_3_3}
\end{figure}

In this numerical example, the goal is to find the right-hand side vector $F_{kh}^{\guess} \in \mathbb{R}^{n_k(n_h-1)}$
with
\begin{align*}
    A_{kh} U_{kh}^{\true} \approx F_{kh}^{\guess}.    
\end{align*}
Therefore, we need to find some initial data $u^0(x)$ and force function $f(x, t)$ that descibes well the observation data $U_{kh}$.
In our loss function, we use Tikhonov regularization with $\alpha = 0.01$ 
and 
\begin{align*}
    J(F_{kh}^{\guess}) := \|U_{kh}(F_{kh}^{true}) - U_{kh}(F_{kh}^{\guess}) \|_2 + \alpha \|F_{kh}^{\guess} \|_2. 
\end{align*}
For the space and time discretization, we choose $n_h = 150$ and $n_k = 50$, which yields a $7450$ dimensional optimization problem.
After $1000$ optimization steps, we reach a loss of $0.0103$, 
whereas the start loss was $60.96$. In Figure \ref{fig:recovered_rhs_1+1d_heat_center} we see that the true solution and the guess of the solution almost coincide for $x = \frac{1}{2}$.

\subsection{1+1D time-stepping heat equation with vector-valued initial condition}
\label{sec:1+1d_time_stepping_heat}
In the fourth numerical test, instead of using a space-time discretization of the 1+1D heat equation, we use a backward Euler time-stepping scheme, which is equivalent to a decoupled solution of equation (\ref{eq:heat_1+1D_space_time_system}).
In contrast to the previous experiment,
the goal of this example is slightly different, as we want to estimate $U_{h}^{0, \guess}$, i.e., the initial condition. 
Once again, the loss function tracks the difference between the observed data $U_{kh}^{\guess}$ and the true (observed) data $U_{kh}^{\true}$.
We note that by using the same optimization formulation as in the previous example the final result would be numerically equivalent. This is because the backward Euler scheme can be derived from a sequential decoupling of the previous space-time discretization.

\begin{figure}[H]
    \centering
    \subfloat[Full]{
        \includegraphics[width=8cm]{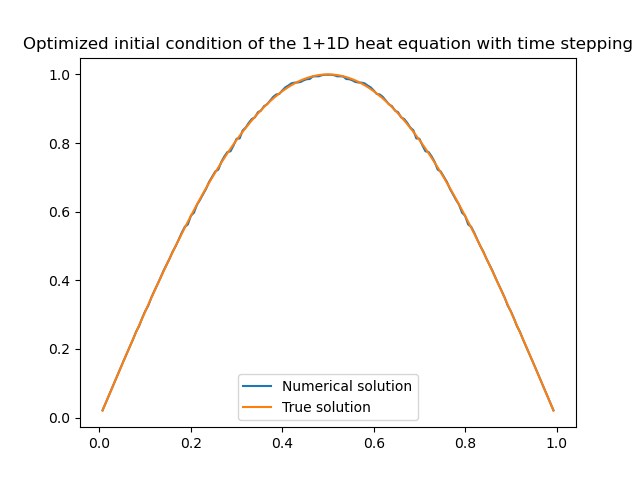}
    }%
    \quad
    \subfloat[Zoomed in]{
        \includegraphics[width=8cm]{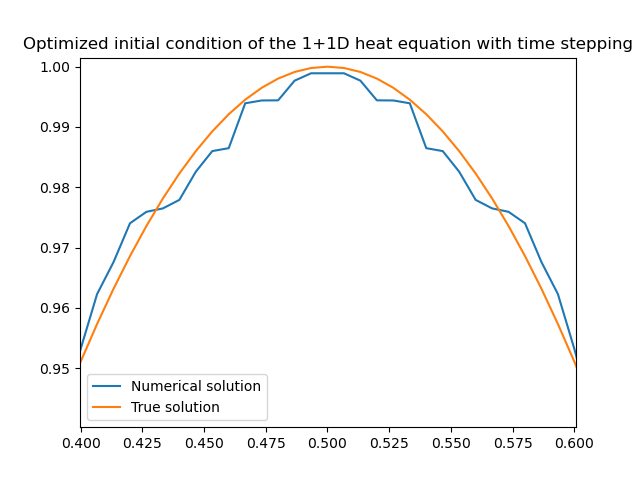}
    }
    \caption{Example 4: Recovered initial condition for the time-stepping formulation of the 1+1D heat equation problem.}\label{fig:recovered_ic_1+1d_heat}
\end{figure}

\begin{figure}[H]
    \hspace{-2cm}
    \includegraphics[width=21cm]{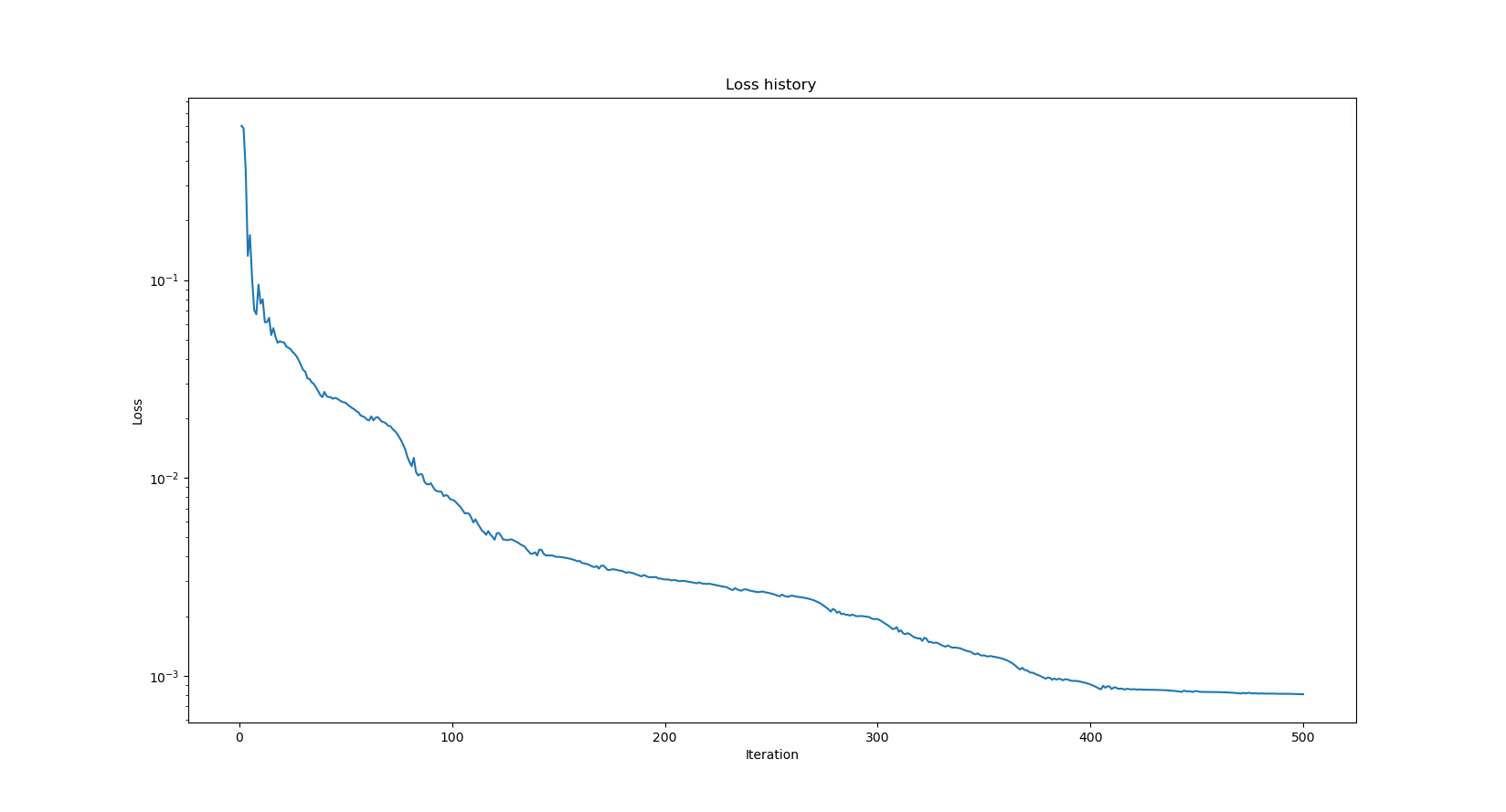}
    \caption{Example 4: Loss history of 1+1D time-stepping heat equation with vector-valued initial condition.}\label{fig:loss_3_4}
\end{figure}

The time-stepping scheme reads as
\begin{align*}
    \left[\frac{1}{k}I_h + K_h\right]U_h^j = F_h^j + \frac{1}{k}U_h^{j-1} \qquad \forall 1 \leq j \leq n_k,
\end{align*}
where the time index is denoted as a superscript
for better readability.
In this example, we decided to infer the exact initial condition $U_h^0$ only and use the loss function
\begin{align*}
     J(U_{h}^{0, \guess}) := \frac{1}{n_k+1}\sum_{j= 0}^{n_k}\|U_{kh}^{j, \true}(U_{h}^{0, \true}) - U_{kh}^{j, \guess}(U_{h}^{0, \guess}) \|_2 + \alpha \|U_{kh}^{0,\guess} \|_2
\end{align*}
with $\alpha = 0.1$.
We use the same temporal and spatial mesh as in Section~\ref{sec:1+1d_space_time_heat}; however, in this case, 
we only optimize the parameters of the initial condition, resulting in $149$ tunable parameters.
After $500$ optimization steps, we reach a loss of $8.05 \cdot 10^{-4}$, see Figure \ref{fig:loss_3_4}. Moreover, we obtain the initial condition shown in Figure \ref{fig:recovered_ic_1+1d_heat}.

\subsection{2D Poisson problem for thermal fin with subdomain-dependent heat conductivities}
\label{sec:2d_thermal_fin}
For the fifth numerical test, we consider a two-dimensional Poisson problem describing heat dissipation in a thermal fin \cite{Qian2017}[Sec. 5.1].
The goal is a classical parameter estimation problem. From given 
reference data, multiple unknown parameters shall be approximated.

The strong form reads: Find $u: \Omega \rightarrow \mathbb{R}$ such that
\begin{align*}
    \sum_{i = 0}^4 - \nabla \cdot (\kappa_i 1_{\Omega_i}(x) \nabla u(x)) &= 0, \qquad \forall x \in \Omega, \\
    \kappa_i \nabla u \cdot n + \operatorname{Bi} u &= 0, \qquad \forall x \in \Gamma_R \cap \Omega_i, \quad 0 \leq i \leq 4, \\
    \kappa_0 \nabla u \cdot n &= 1. \qquad \forall x \in \Gamma_N.
\end{align*}
The domain, subdomains and boundaries are given as
\begin{align*}
    \Omega_0 &= (2.5, 3.5) \times (0,4), \\
    \Omega_i &= ((0, 2.5) \cup (3.5, 6)) \times (i-0.25, i), \qquad \forall 1 \leq i \leq 4, \\
    \Omega &= \bigcup_{i =
    0}^4 \Omega_i, \\
    \Gamma_N &= (2.5, 3.5) \times \{ 0 \}, \\
    \Gamma_R &= \partial \Omega \setminus \Gamma_N.
\end{align*}

\begin{figure}[H]
    \centering
    \subfloat[Mesh]{
        \includegraphics[width=7cm]{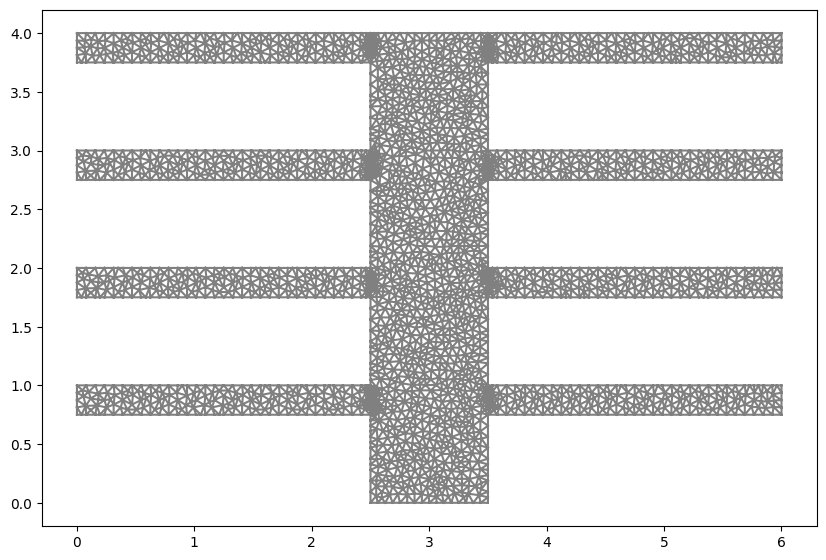}
    }%
    \quad
    \subfloat[Subdomains]{
        \includegraphics[width=7cm]{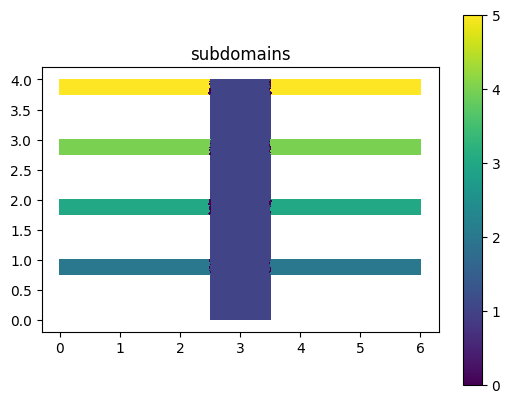}
    }
    \caption{Example 5: Computational domain of 2D Poisson problem for thermal fin.}\label{fig:thermal_fin_domain}
\end{figure}

\begin{figure}[H]
    \centering
    \subfloat[True Solution]{
        \includegraphics[width=8cm]{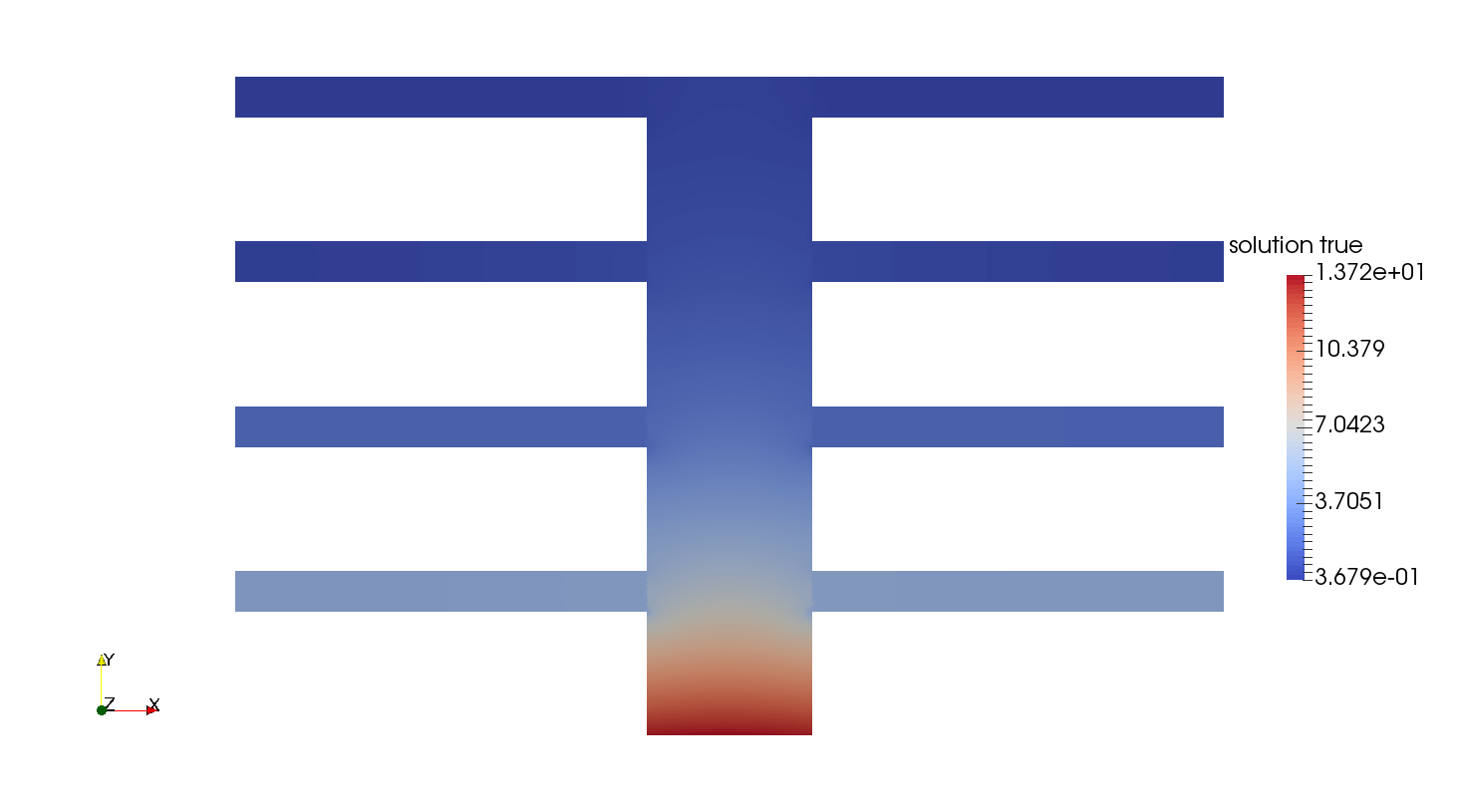}
    }%
    \quad
    \subfloat[Recovered Solution]{
        \includegraphics[width=8cm]{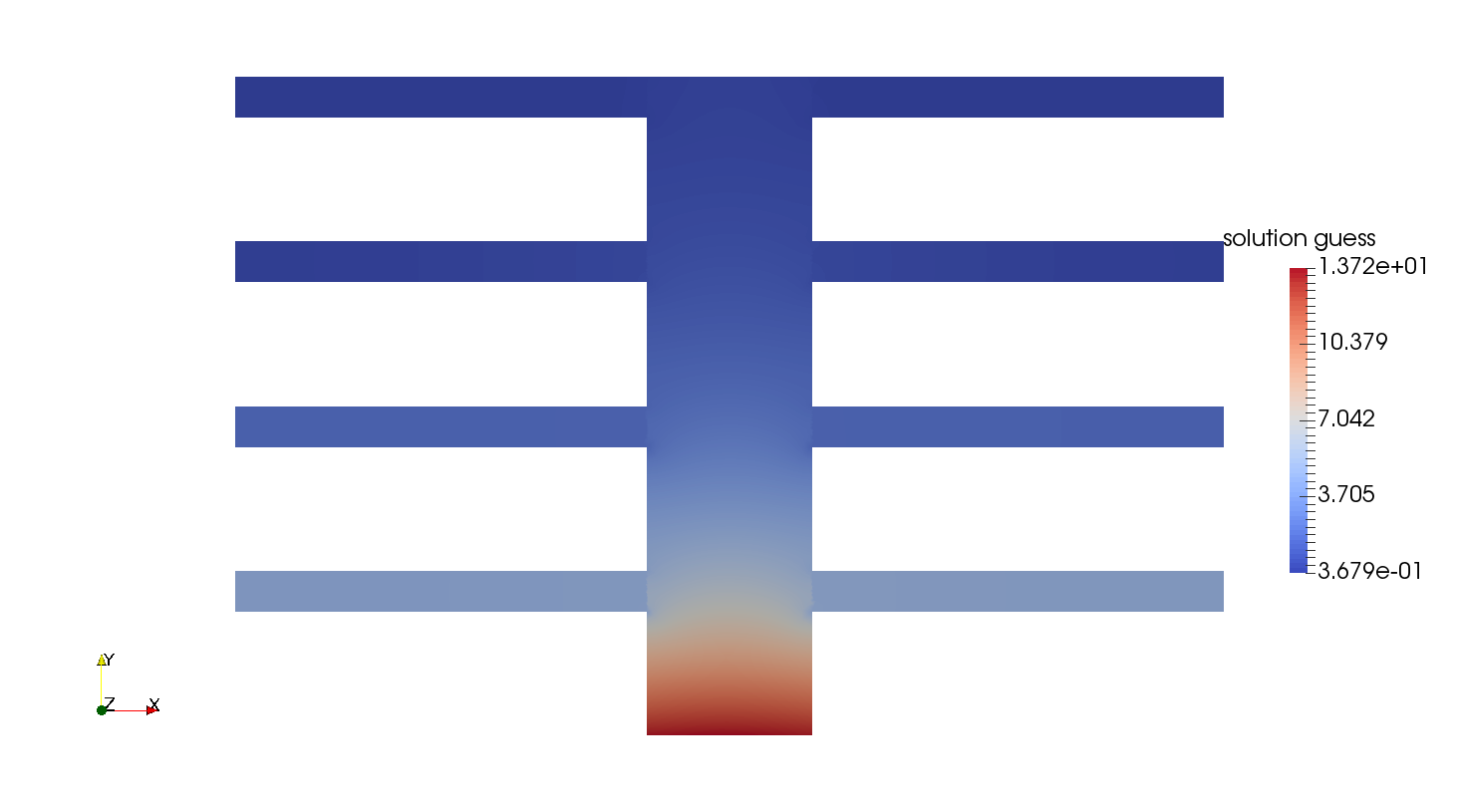}
    }
    \quad
    \subfloat[Solution Difference]{
        \includegraphics[width=8cm]{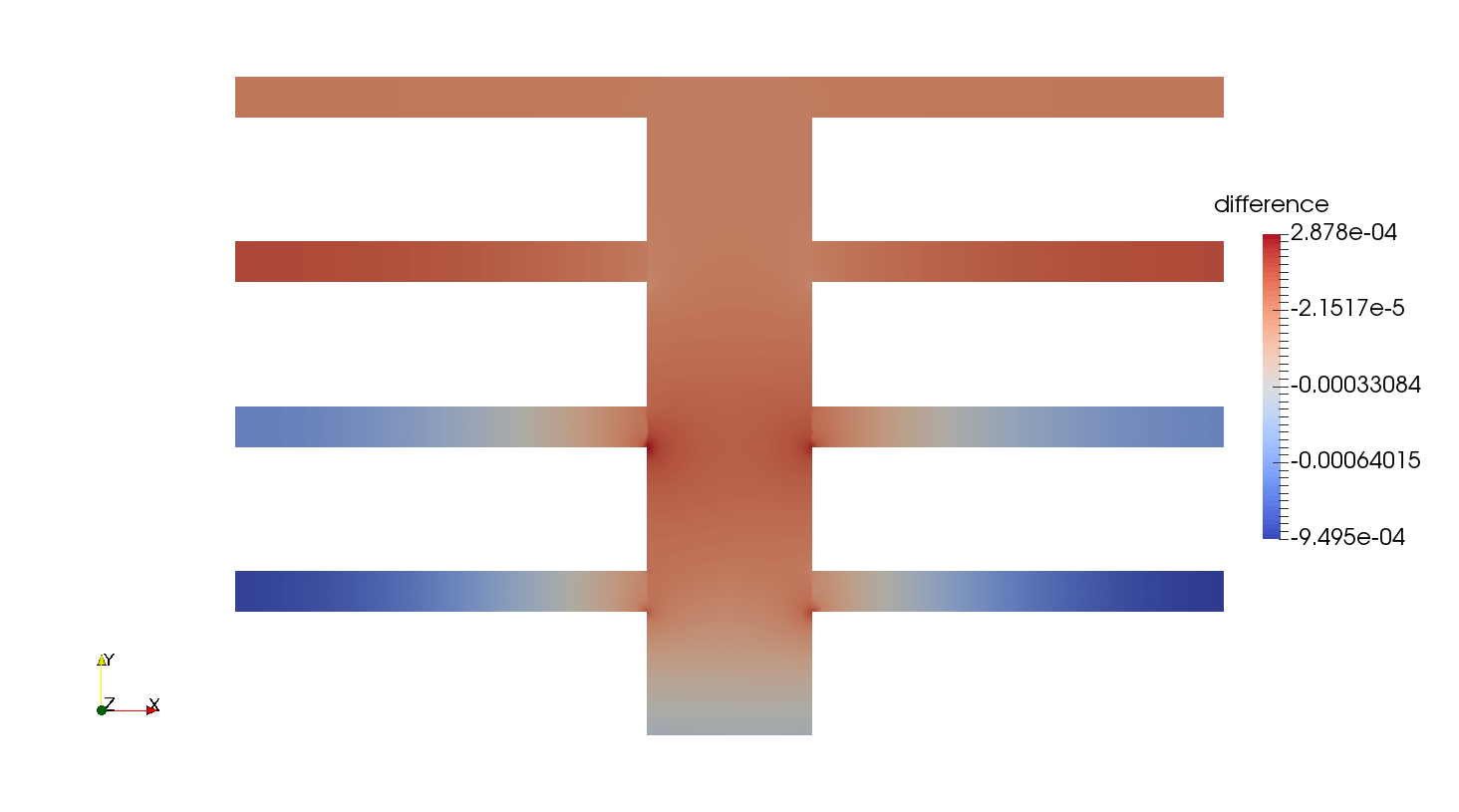}
    }
    \quad
    \subfloat[Loss history]{
        \includegraphics[width=8cm]{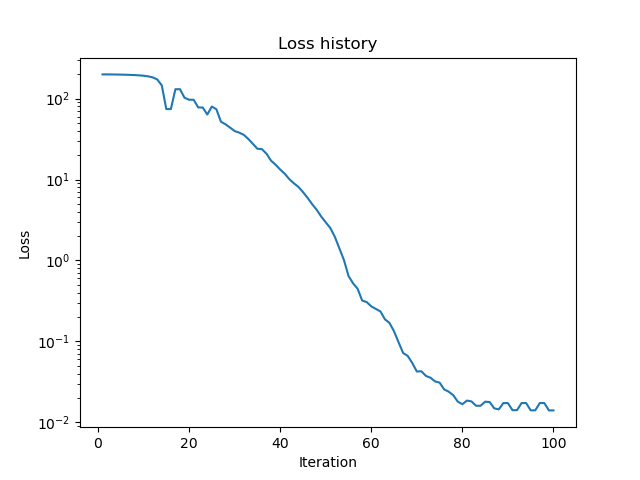}
    }
    \caption{Example 5: Solutions of 2D Poisson problem for thermal fin after optimization and Loss-history.}
    \label{fig:thermal_fin_solution}
    \label{fig:loss_3_5}
\end{figure}

The parameters that need to be learned are given as
\begin{align*}
    \mu = \begin{pmatrix}
        \kappa_0, \kappa_1, \kappa_2, \kappa_3, \kappa_4, \operatorname{Bi} 
    \end{pmatrix}^T \in \mathbb{R}^6 \qquad \text{with} \qquad \kappa_i \in [0.1, 10], \forall 0 \leq i \leq 4, \quad \operatorname{Bi} \in [0.01, 1].
\end{align*}
The weak form for this problem reads: Find $u \in H^1(\Omega)$ such that
\begin{align*}
    \sum_{i = 0}^4 \kappa_i (\nabla u, \nabla \phi)_{\Omega_i} + \operatorname{Bi} (u, \phi)_{\Gamma_R} = (1, \phi)_{\Gamma_N} \qquad \forall \phi \in H^1(\Omega).
\end{align*}
For the experiment, we use the true solution
\begin{align*}
    \mu_{\true} = \begin{pmatrix}
        0.1, 8.37317, 6.57228, 0.466517, 1.88354, 0.01
    \end{pmatrix}^T    
\end{align*}
and the reference solution
\begin{align*}
    \mu_{\operatorname{ref}} = \begin{pmatrix}
        1, 1, 1, 1, 1, 0.1
    \end{pmatrix}^T,
\end{align*}
which were both taken from \cite{Keil2021}.
As an initial guess, we use
\begin{align*}
    \mu_{\guess} = \begin{pmatrix}
        0.5, 0.5, 0.5, 0.5, 0.5, 0.5
    \end{pmatrix}^T.
\end{align*}
We define the loss function
\begin{align*}
    J(\mu^{\guess}) := \|U_h(\mu^{\true}) - U_h(\mu^{\guess})\|_2 + 0.1 \left\|\frac{\mu^{\guess}-\mu^{\operatorname{ref}}}{\mu^{\operatorname{ref}}}\right\|_2
\end{align*}
and clip the parameter guess $\mu^{\guess}$ to the parameter space $[0.1, 10]^5 \times [0.01, 1]$.
The linear finite element space consists of $3344$ degrees of freedom.
After $100$ optimization steps with Rprop and a learning rate of $0.01$, we achieve a loss of $0.014$, see Figure \ref{fig:loss_3_5} (d), and get the parameter guess
\begin{align*}
    \mu^{\guess} = \begin{pmatrix}
        0.1, 8.3231, 6.5108, 0.4669, 1.8841, 0.01
    \end{pmatrix}^T.
\end{align*}

\subsection{2+1D nonlinear heat equation with unknown initial condition}
\begin{figure}[H]
    \centering
    \subfloat[True Initial Condition]{
        \includegraphics[width=8cm]{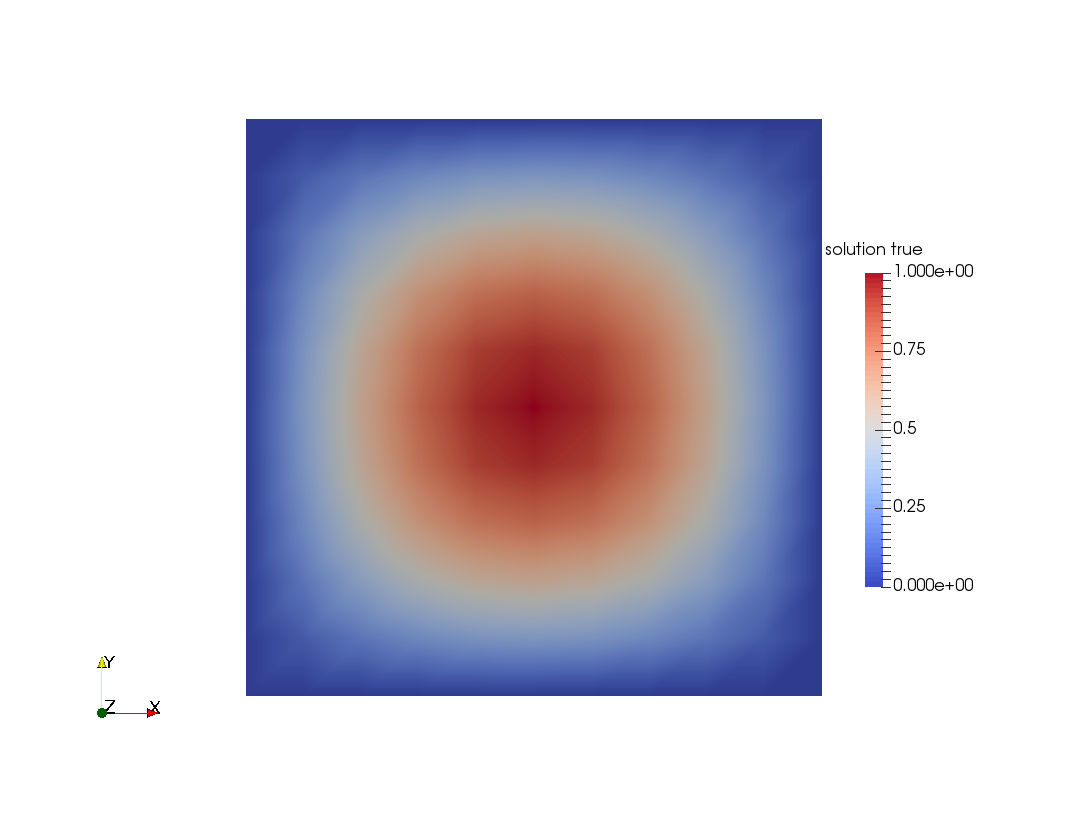}
    }%
    \quad
    \subfloat[Recovered Initial Condition]{
        \includegraphics[width=8cm]{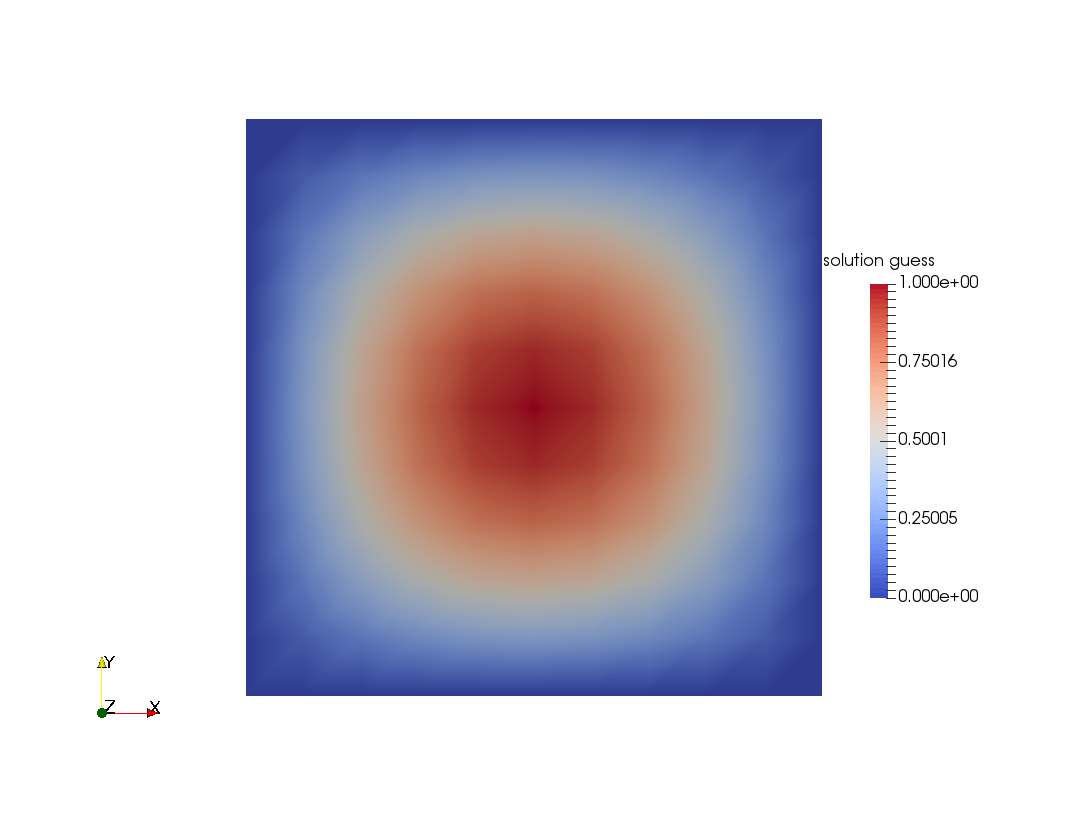}
    }
    \quad
    \subfloat[Initial Condition Difference]{
        \includegraphics[width=8cm]{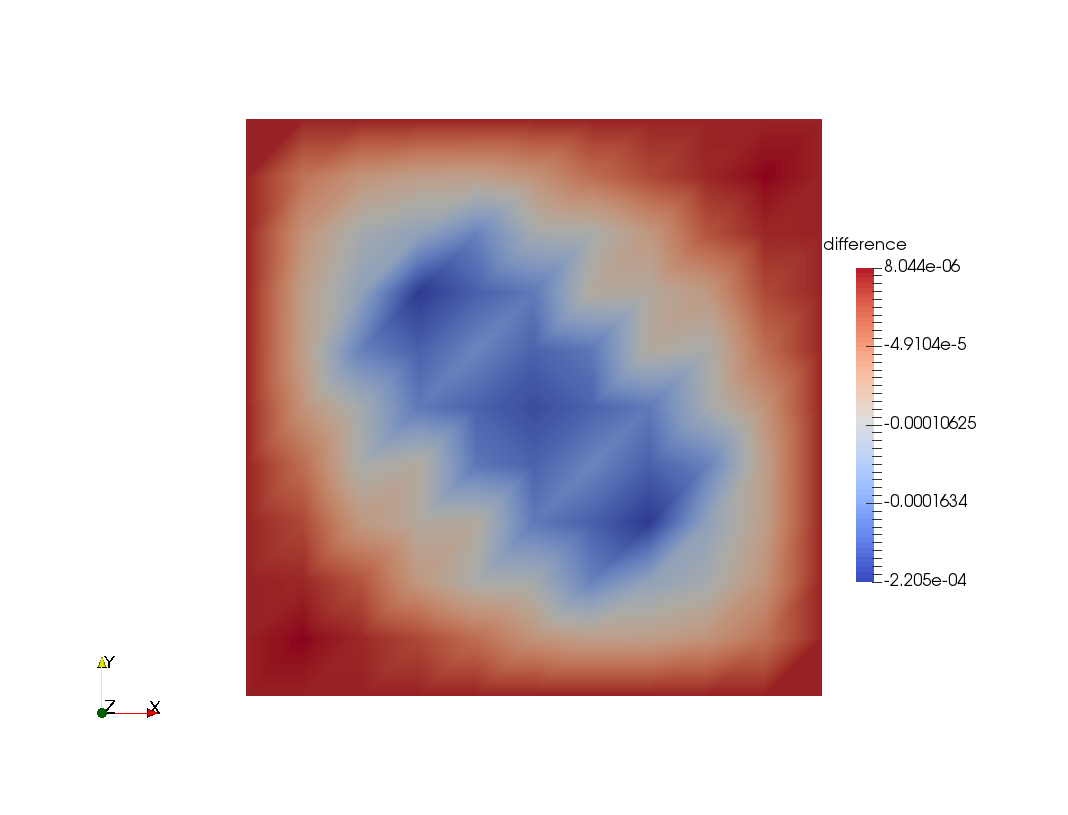}
    }
    \quad
    \subfloat[Recovered solution]{
    \includegraphics[width=8cm]{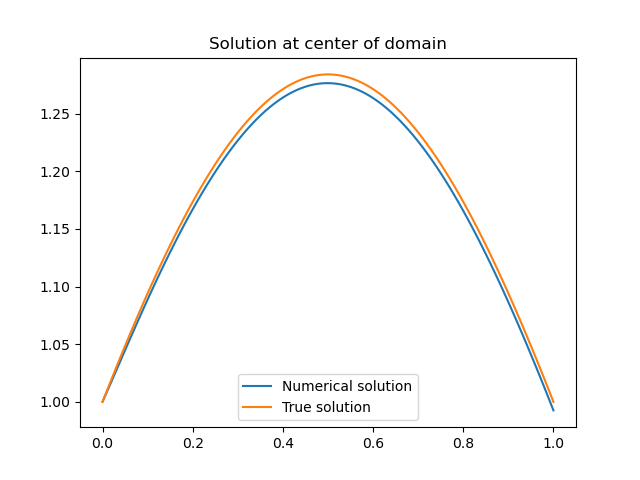}
    }
    \caption{Example 6: Initial conditions of 2+1D nonlinear heat equation after optimization (a)--(c). 
    Recovered solution of the nonlinear heat equation in 2+1D at $(x,y) = (\frac{1}{2},\frac{1}{2})$ (d).}
    \label{fig:nonlinear_heat_solution}
    \label{fig:recovered_rhs_2+1d_nonlinear_heat_center}
\end{figure}

\begin{figure}[H]
    \centering
    \includegraphics[width=14.5cm]{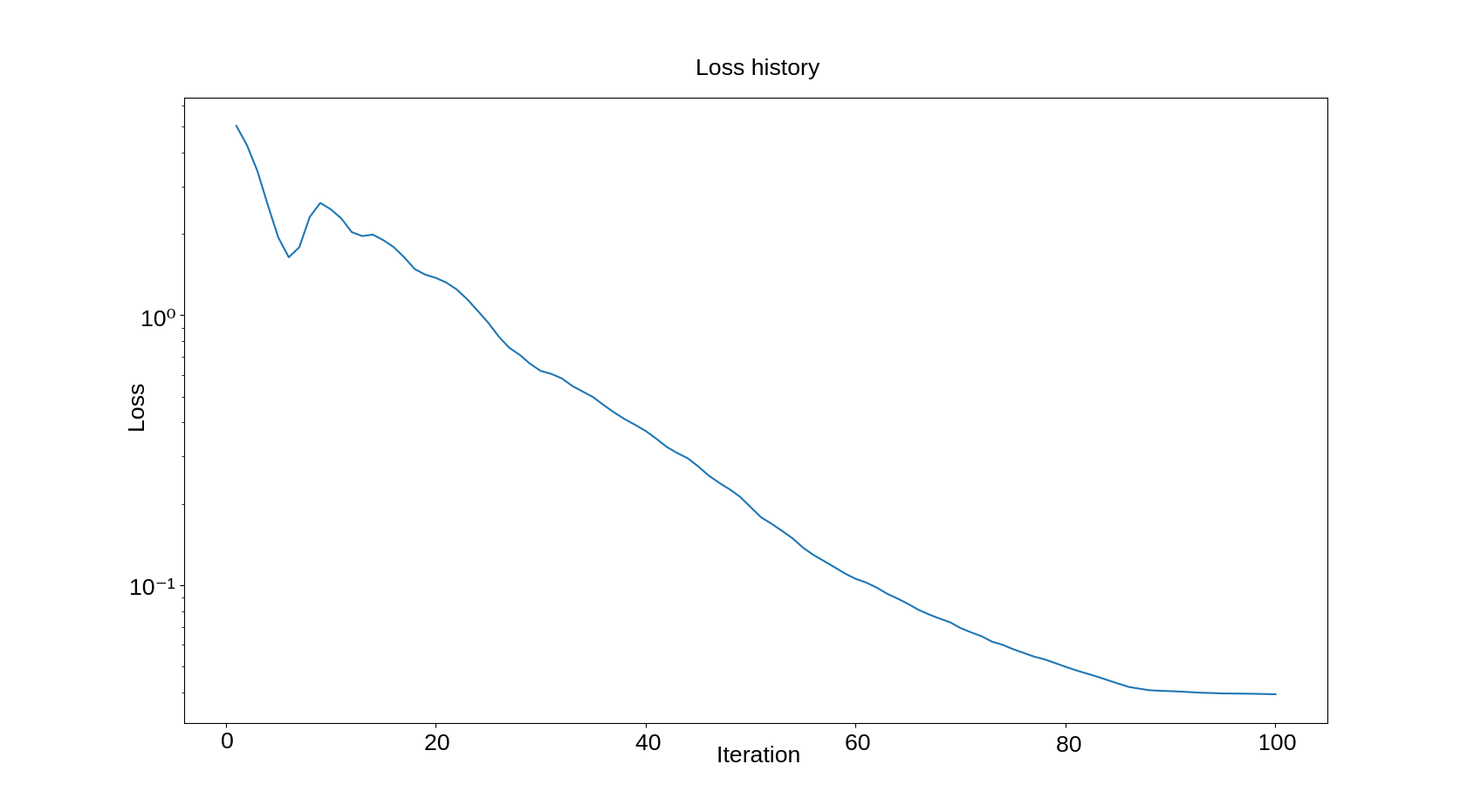}
    \caption{Example 6: Loss history of 2+1D nonlinear heat equation with unknown initial condition.}
    \label{fig:loss_3_6}
\end{figure}

On our sixth experiment, we consider a modification of the nonlinear heat equation example from DOpElib \cite{DOpElib}. 
The goal is to learn the initial condition $u^{0,\true}$ from given 
observations compared to reference data.
The strong form reads: Find $u: \Omega \times I \rightarrow \mathbb{R}$ such that
\begin{align*}
    \partial_t u  - \Delta u + u^2 &= f, \qquad \text{in }  \Omega \times I, \\
    u &= 0, \qquad \text{on } \partial \Omega \times I, \\
    u &= u^0 \qquad \text{on } \Omega \times \{ 0 \}.
\end{align*}
The space-time domain is given by
\begin{align*}
    \Omega &:= (0,1)^2, \\
    I &:= (0,T) := (0,1).
\end{align*}
After the application of a time discretization, i.e., the Crank-Nicolson method with the uniform time step size $k := t_n - t_{n-1}$, the weak form for this problem reads: Find $u_n \in H^1_0(\Omega) := \{u \in H^1(\Omega) | u = 0 \text{ on } \partial \Omega \}$ such that
\begin{align*}
    (u_n, \phi) + \frac{k}{2}(\nabla u_n, \nabla \phi) + \frac{k}{2}(u_n^2, \phi) = (u_{n-1}, \phi) - \frac{k}{2}(\nabla u_{n-1}, \nabla \phi) - \frac{k}{2}(u_{n-1}^2, \phi) + \frac{k}{2}(f_{n-1} + f_n, \phi) \\
    \forall \phi \in H^1_0(\Omega).
\end{align*}
Using the force function
\begin{align*}
    f(x, y, t) := \left(- 2 t e^{t^{2}} + e^{t} \sin{\left(\pi x \right)} \sin{\left(\pi y \right)} + e^{t^{2}} + 2 \pi^{2} e^{t^{2}}\right) e^{- 2 t^{2} + t} \sin{\left(\pi x \right)} \sin{\left(\pi y \right)}
\end{align*}
and the space-time solution
\begin{align*}
    u^{\true}(x, y, t) := \exp(t-t^2)\sin(\pi x)\sin(\pi y),
\end{align*}
we want to learn the initial condition which is given by
\begin{align*}
    u^{0,\true}(x, y) = \sin(\pi x)\sin(\pi y).
\end{align*}
As a first guess for the initial condition, we use
\begin{align*}
    u^{0,\guess}(x, y) = 0.
\end{align*}
Finally, we define the loss function
\begin{align*}
    J(u^{0,\guess}) := \|U_h^{0,\true} - U_h^{0,\guess}\|_2 + \|U_h^{T,\true} - U_h^{T,\guess}(U_h^{0,\guess})\|_2 + 0.1 \left\| U_h^{0,\guess} \right\|_2.
\end{align*}
For the spatial discretization, we use linear finite elements with 10 elements in x- and y-direction, leading to 121 tunable parameters.
In time, we have 100 time steps.
After 100 optimization steps with Rprop and a learning rate of $0.1$, we achieve a loss of $0.0393$. In Figure \ref{fig:recovered_rhs_2+1d_nonlinear_heat_center}, we plot the solution at the center of the domain over time. 
The decrease of the cost functional over the optimization steps is presented in Figure \ref{fig:loss_3_6}.

\subsection{2D Navier-Stokes with boundary control to minimize drag}
In the seventh experiment, we use the stationary Navier-Stokes example from \cite{Becker2001OC}, which is a modification of the Navier-Stokes 2D-1 benchmark \cite{SchaeferTurek1996} with Neumann boundary control on parts of the walls above and below the cylindrical obstacle.
The strong form reads: 
\begin{align*}
    -\nabla \cdot \sigma + (v \cdot \nabla)v &= 0 \qquad \text{in }  \Omega, \\
    \nabla \cdot v &= 0 \qquad \text{in }  \Omega,
\end{align*}
with the unsymmetric stress tensor 
\begin{align*}
    \sigma := \sigma(v, p) := -pI + \nu \nabla v.   
\end{align*}
Consequently, we have to find the vector-valued velocity $v: \Omega \rightarrow \mathbb{R}^2$ and the scalar-valued pressure $p: \Omega \rightarrow \mathbb{R}$ such that
\begin{align*}
    - \nu \Delta v + \nabla p + (v \cdot \nabla)v  &= 0 \qquad \text{in }  \Omega, \\
    \nabla \cdot v &= 0 \qquad \text{in }  \Omega,
\end{align*}
where, $\nu = 10^{-3}$ denotes the kinematic viscosity. 
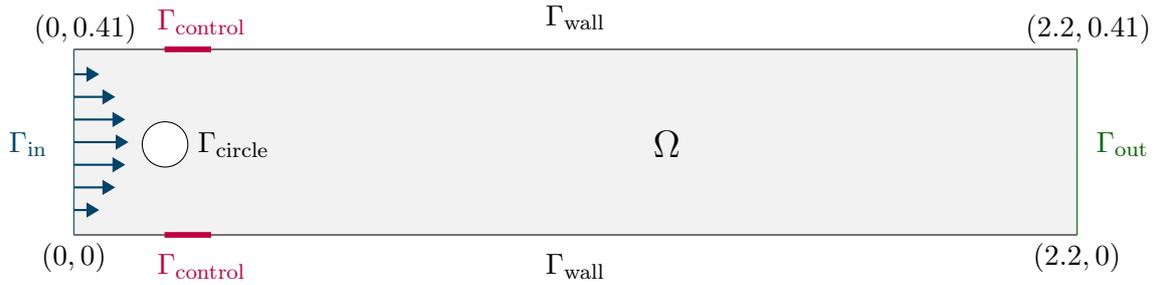
\begin{figure}[H]
    \begin{center}
    \begin{tikzpicture}[scale = 6, draw=black]
        \draw[draw=white, fill=black!5!white] (0,0) rectangle (2.2,0.41);
        \draw[draw=black] (0,0) -- (2.2,0);
        \draw[draw=black] (0,0.41) -- (2.2,0.41);
        \draw[draw=blue!60!black] (0,0) -- (0,0.41);
        \draw[draw=green!40!black] (2.2,0) -- (2.2,0.41);
        
        \draw[fill=white] (0.2,0.2) circle (0.05);
        \node (circ) at (0.35,0.2) {$\Gamma_{\text{circle}}$};
        
        \node (omega) at (1.3,0.2) {\Large{$\Omega$}};
        
        \node (in) at (-0.1,0.2) {{\color{blue!60!black}$\Gamma_{\text{in}}$}};
        \node (out) at (2.3,0.2) {{\color{green!40!black}$\Gamma_{\text{out}}$}};
        \node (wall-one) at (1.1,0.48) {$\Gamma_{\text{wall}}$};
        \node (wall-two) at (1.1,-0.07) {$\Gamma_{\text{wall}}$};
        
        \node (corner-one) at (0,-0.05) {$(0,0)$};
        \node (corner-two) at (0.025,0.455) {$(0,0.41)$};
        \node (corner-three) at (2.2,-0.05) {$(2.2,0)$};
        \node (corner-four) at (2.225,0.455) {$(2.2,0.41)$};
        
        \draw[draw=blue!60!black, -Triangle, thick] (0,0.055) -- (0.056,0.055);
        \draw[draw=blue!60!black, -Triangle, thick] (0,0.105) -- (0.091,0.105);
        \draw[draw=blue!60!black, -Triangle, thick] (0,0.155) -- (0.113,0.155);
        \draw[draw=blue!60!black, -Triangle, thick] (0,0.205) -- (0.12,0.205);
        \draw[draw=blue!60!black, -Triangle, thick] (0,0.255) -- (0.113,0.255);
        \draw[draw=blue!60!black, -Triangle, thick] (0,0.305) -- (0.091,0.305);
        \draw[draw=blue!60!black, -Triangle, thick] (0,0.355) -- (0.056,0.355);

        \draw[draw=purple, fill=purple] (0.2,-0.005) rectangle (0.3,0.005);
        \draw[draw=purple, fill=purple] (0.2,0.405) rectangle (0.3,0.415);
        \node (control-one) at (0.28, 0.47) {{\color{purple}$\Gamma_{\text{control}}$}};
        \node (control-two) at (0.28, -0.07) {{\color{purple}$\Gamma_{\text{control}}$}};

    \end{tikzpicture}
    \caption{Example 7: Domain of the Navier-Stokes benchmark problem.}
    \label{fig:channel_domain}
    \end{center}
\end{figure}

\begin{figure}[H]
    \centering
    \subfloat[Loss history]{
        \includegraphics[width=8cm]{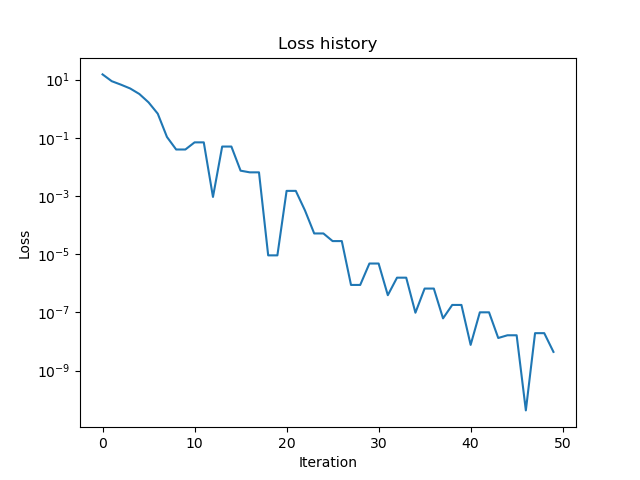}
    }
    \subfloat[Boundary control]{
        \includegraphics[width=8cm]{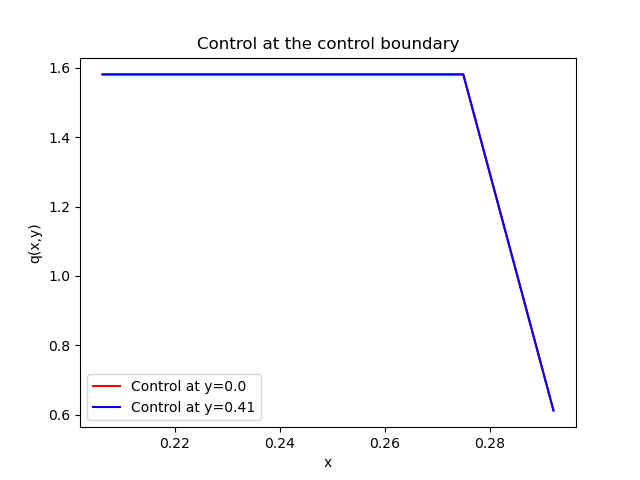}
    }
    \caption{Example 7: Loss history of Neumann boundary control for stationary Navier-Stokes equations (a) and the optimal
    Neumann boundary control $q^{\guess}$ for Navier-Stokes after 50 optimization steps.}
    \label{fig:loss_navier_stokes}
    \label{fig:optimal_control_navier_stokes}
\end{figure}

The domain is defined as $\Omega := (0,2.2) \times (0,0.41) \setminus B_r(0.2,0.2)$ with $r = 0.05$, where $B_r(x,y)$ denotes the ball around $(x,y)$ with radius $r$. On the boundary $\partial \Omega$ we prescribe Dirichlet and Neumann boundary conditions. We apply the no-slip condition, i.e. homogeneous Dirichlet boundary conditions $v = 0$ on $\Gamma_{\text{wall}} \cup \Gamma_{\text{circle}}$, where ${\Gamma_{\text{wall}} = \left((0,0.2) \cup (0.3, 2.2)\right) \times \lbrace 0, 0.4.1 \rbrace}$ and ${\Gamma_{\text{circle}} = \partial B_r(0.2,0.2)}$. On the boundary ${\Gamma_{\text{in}} = \lbrace 0 \rbrace \times (0,0.41)}$, we enforce a parabolic inflow profile, which is a type of inhomogeneous Dirichlet boundary condition, i.e. $v = g$ on $\Gamma_{\text{in}}$. The inflow parabola is given by
\begin{align*}
	g(0,y) = \begin{pmatrix}
		\frac{1.2y(0.41-y)}{0.41^2} \\
		0
	\end{pmatrix},
\end{align*}
For brevity, we denote the Dirichlet boundary as $\Gamma_D := \Gamma_{\text{in}} \cup \Gamma_{\text{wall}} \cup \Gamma_{\text{circle}}$.
The outflow boundary $\Gamma_{\text{out}} = \lbrace 2.2 \rbrace \times (0,0.41)$ has Neumann boundary conditions, i.e. $\sigma \cdot n = 0$ on $\Gamma_{\text{out}}$.
The goal of our optimization loop is to control the inhomogeneous Neumann boundary condition on $\Gamma_{\text{control}} = (0.2, 0.3) \times \lbrace  0, 0.41 \rbrace$ such that the drag coefficient is minimized.
 The drag coefficient \cite{SchaeferTurek1996} is defined as
\begin{align*}
    C_D(v, p) = 500 \int_{\Gamma_{\text{circle}}}\sigma(v, p) \cdot n \cdot \begin{pmatrix}
        1 \\ 0
    \end{pmatrix}\ \mathrm{d}x.
\end{align*}
For the Neumann boundary condition, we need to find a function $q: \Gamma_{\text{control}} \rightarrow \mathbb{R}$ such that $\sigma \cdot n = qn$ on $\Gamma_{\text{control}}$.
For the implementation of the weak formulation in \texttt{FEniCS}, we refer to the old \texttt{FEniCS} tutorials\footnote{\url{https://fenicsproject.org/olddocs/dolfin/1.6.0/python/demo/documented/stokes-iterative/python/documentation.html}}.
In short, the weak form is given by:
Find $v \in \{g + H^1_{0, \Gamma_D}(\Omega)^2\}$ and $p \in L^2(\Omega)$ such that for all test functions
$\phi^v \in H^1_{0, \Gamma_D}(\Omega)^2 := \{f \in H^1(\Omega)^2 | f = 0 \text{ on } \Gamma_D\}$ and $ \phi^p \in L^2(\Omega)$ it holds
\begin{align*}
    \nu (\nabla v, \nabla \phi^v) - (p, \nabla \phi^v) + ((v \cdot \nabla)v, \phi^v) + (\nabla \cdot v, \phi^p) = \langle q, \phi^v \cdot n\rangle_{\Gamma_{\text{control}}}.
\end{align*}
As a first guess for the control, we use $q^{\guess}(x, y) = 0$ and enforce $q^{\guess} \geq 0$, as well as $q^{\guess}(0.3, y) \leq 0.7$ (to avoid divergence of Newton solver) during the optimization process.
We use the quadratic loss function
\begin{align*}
    J(q^{\guess}) := C_D(v(q^{\guess}), p(q^{\guess}))^2.
\end{align*}
For the spatial discretization, we use inf-sup stable Taylor-Hood elements with quadratic elements for velocity and linear elements for pressure with a total of 21400 DoFs. Nevertheless, since we only control parts of the boundary, we have merely 12 tunable parameters. For optimization, we use Rprop with a learning rate of $0.1$ and show the history of the loss over the training in Figure \ref{fig:loss_navier_stokes}. 
After 50 iterations, we end up with the boundary control $q^{\guess}$ as shown in Figure \ref{fig:optimal_control_navier_stokes} together with the loss history.

\begin{figure}[H]
    \centering
    \subfloat[Velocity Magnitude]{
        \includegraphics[width=8cm]{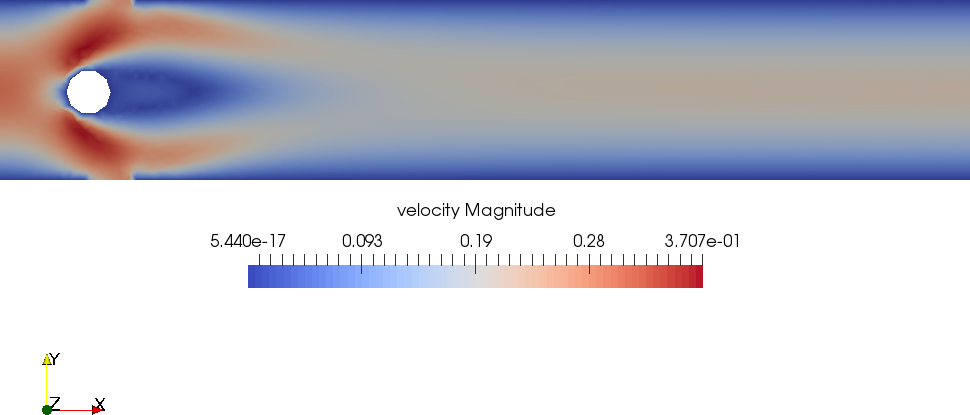}
    }%
    \subfloat[Pressure]{
        \includegraphics[width=8cm]{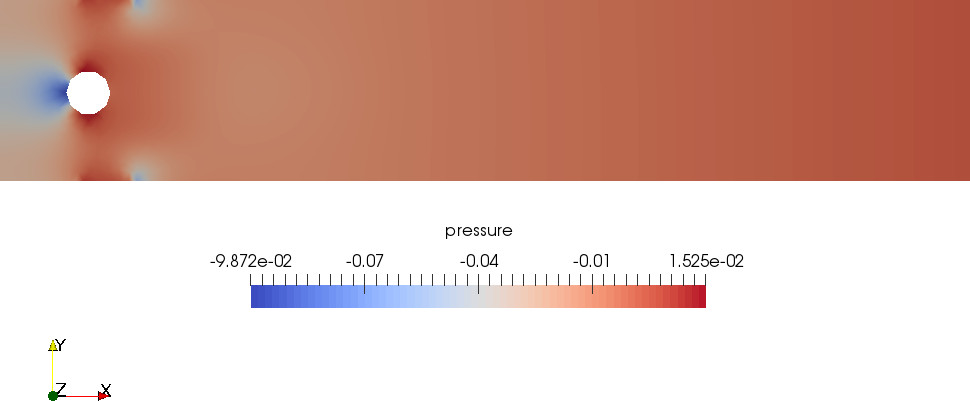}
    }
    \caption{Example 7: Initial Navier-Stokes solution ($q^{\guess} = 0$).}\label{fig:initial_navier_stokes_solution}
\end{figure}
\begin{figure}[H]
    \centering
    \subfloat[Velocity Magnitude]{
        \includegraphics[width=8cm]{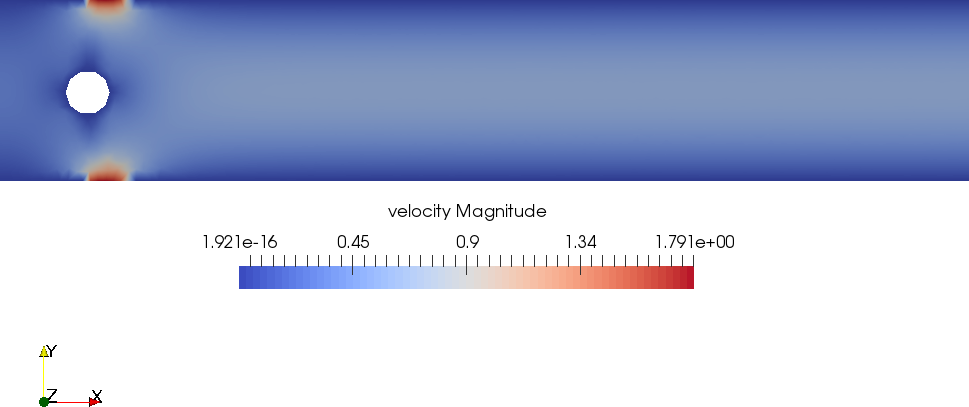}
    }%
    \subfloat[Pressure]{
        \includegraphics[width=8cm]{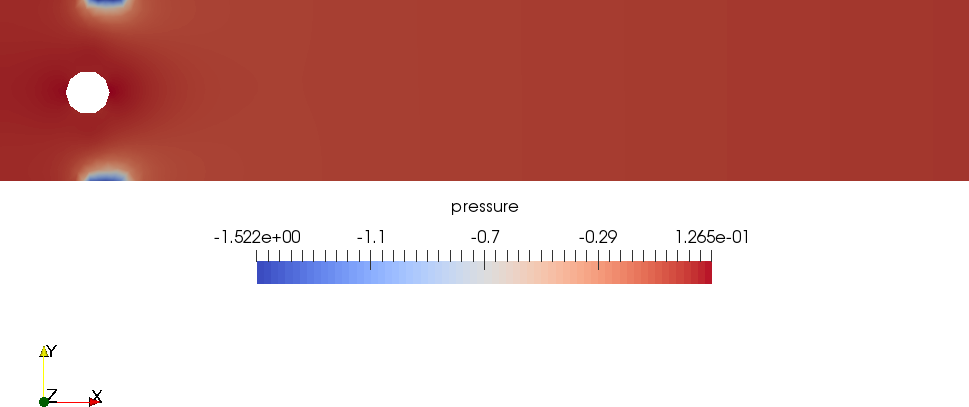}
    }\\
    \subfloat[Velocity x-direction]{
        \includegraphics[width=8cm]{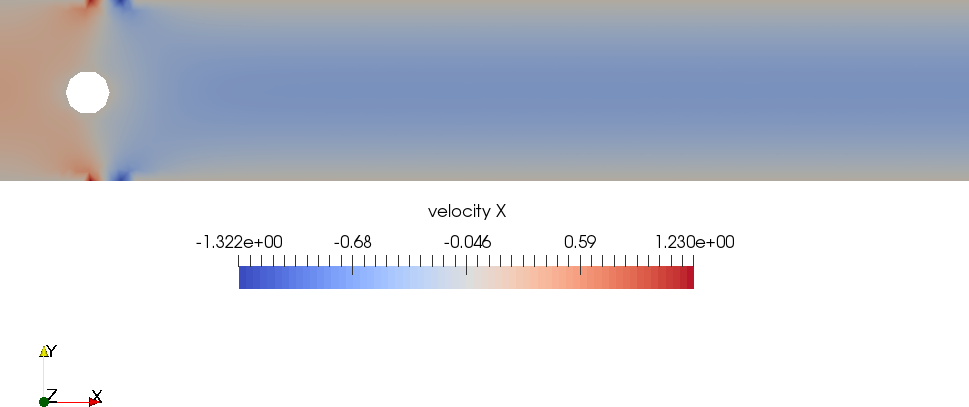}
    }%
    \subfloat[Velocity y-direction]{
        \includegraphics[width=8cm]{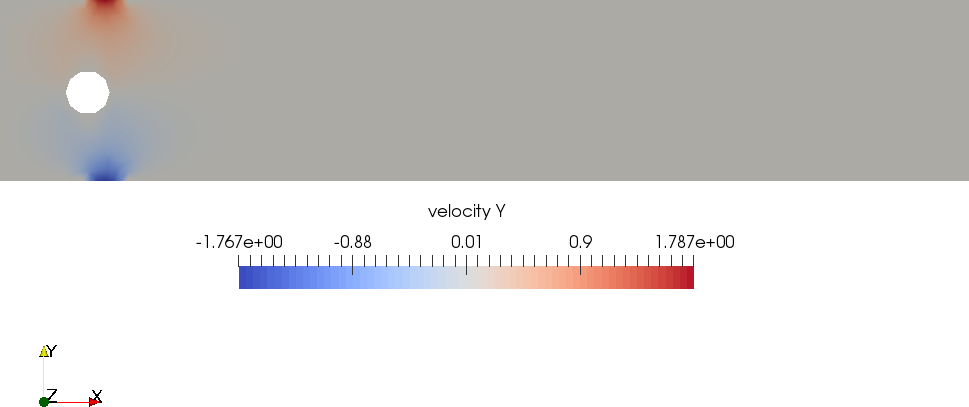}
    }%
    \caption{Example 7: Navier-Stokes solution after 50 optimization steps.}\label{fig:optimized_navier_stokes_solution}
\end{figure}

We observe that our results in Figure \ref{fig:optimized_navier_stokes_solution} coincide with the results obtained by DOpElib \cite{DOpElib} and the fluid is being "sucked out of the domain" at the control boundary $\Gamma_{\text{control}}$.

\subsection{2D Fluid-Structure Interaction with parameter estimation for Lam{\'e} parameter}

\begin{figure}[H]
    \begin{center}
    \begin{tikzpicture}[scale = 6, draw=black]
        \draw[draw=white, fill=black!5!white] (0,0) rectangle (2.2,0.41);
        \draw[draw=black] (0,0) -- (2.2,0);
        \draw[draw=black] (0,0.41) -- (2.2,0.41);
        \draw[draw=blue!60!black] (0,0) -- (0,0.41);
        \draw[draw=green!40!black] (2.2,0) -- (2.2,0.41);
        
        \draw[fill=white] (0.2,0.2) circle (0.05);
        \node (circ) at (0.2,0.1) {$\Gamma_{\text{circle}}$};

        \draw[draw=red, thick] (0.24758, 0.18) -- (0.6, 0.18) -- (0.6, 0.22) -- (0.24758, 0.22);
        \node (beam) at (0.5, 0.13) {{\color{red}$\Gamma_{\text{interface}}$}};

        \draw [draw=orange, thick] (0.24757, 0.18) arc [
    		start angle=-23,
    		end angle=23,
    		x radius=0.05,
    		y radius=0.05
    	];
        \node (beam) at (0.3,0.3) {{\color{orange}$\Gamma_{\text{beam}}$}};

        \draw[fill] (0.6, 0.2) circle [radius=0.006];
        \node (deflect) at (0.725, 0.2) {\scriptsize$(0.6, 0.2)$};
        
        \node (omega) at (1.3,0.2) {\Large{$\Omega$}};
        
        \node (in) at (-0.1,0.2) {{\color{blue!60!black}$\Gamma_{\text{in}}$}};
        \node (out) at (2.3,0.2) {{\color{green!40!black}$\Gamma_{\text{out}}$}};
        \node (wall-one) at (1.1,0.48) {$\Gamma_{\text{wall}}$};
        \node (wall-two) at (1.1,-0.07) {$\Gamma_{\text{wall}}$};
        
        \node (corner-one) at (0,-0.05) {$(0,0)$};
        \node (corner-two) at (0.025,0.455) {$(0,0.41)$};
        \node (corner-three) at (2.2,-0.05) {$(2.5,0)$};
        \node (corner-four) at (2.225,0.455) {$(2.5,0.41)$};
        
        \draw[draw=blue!60!black, -Triangle, thick] (0,0.055) -- (0.056,0.055);
        \draw[draw=blue!60!black, -Triangle, thick] (0,0.105) -- (0.091,0.105);
        \draw[draw=blue!60!black, -Triangle, thick] (0,0.155) -- (0.113,0.155);
        \draw[draw=blue!60!black, -Triangle, thick] (0,0.205) -- (0.12,0.205);
        \draw[draw=blue!60!black, -Triangle, thick] (0,0.255) -- (0.113,0.255);
        \draw[draw=blue!60!black, -Triangle, thick] (0,0.305) -- (0.091,0.305);
        \draw[draw=blue!60!black, -Triangle, thick] (0,0.355) -- (0.056,0.355);
    \end{tikzpicture}
    \caption{Example 8: Domain of the fluid-structure interaction example.}
    \label{fig:fsi_domain}
    \end{center}
\end{figure}
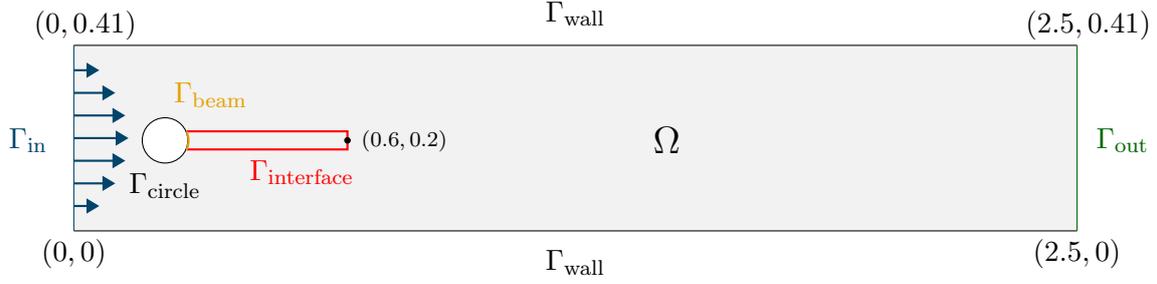

\begin{figure}[H]
    \centering
    \includegraphics[width=12cm]{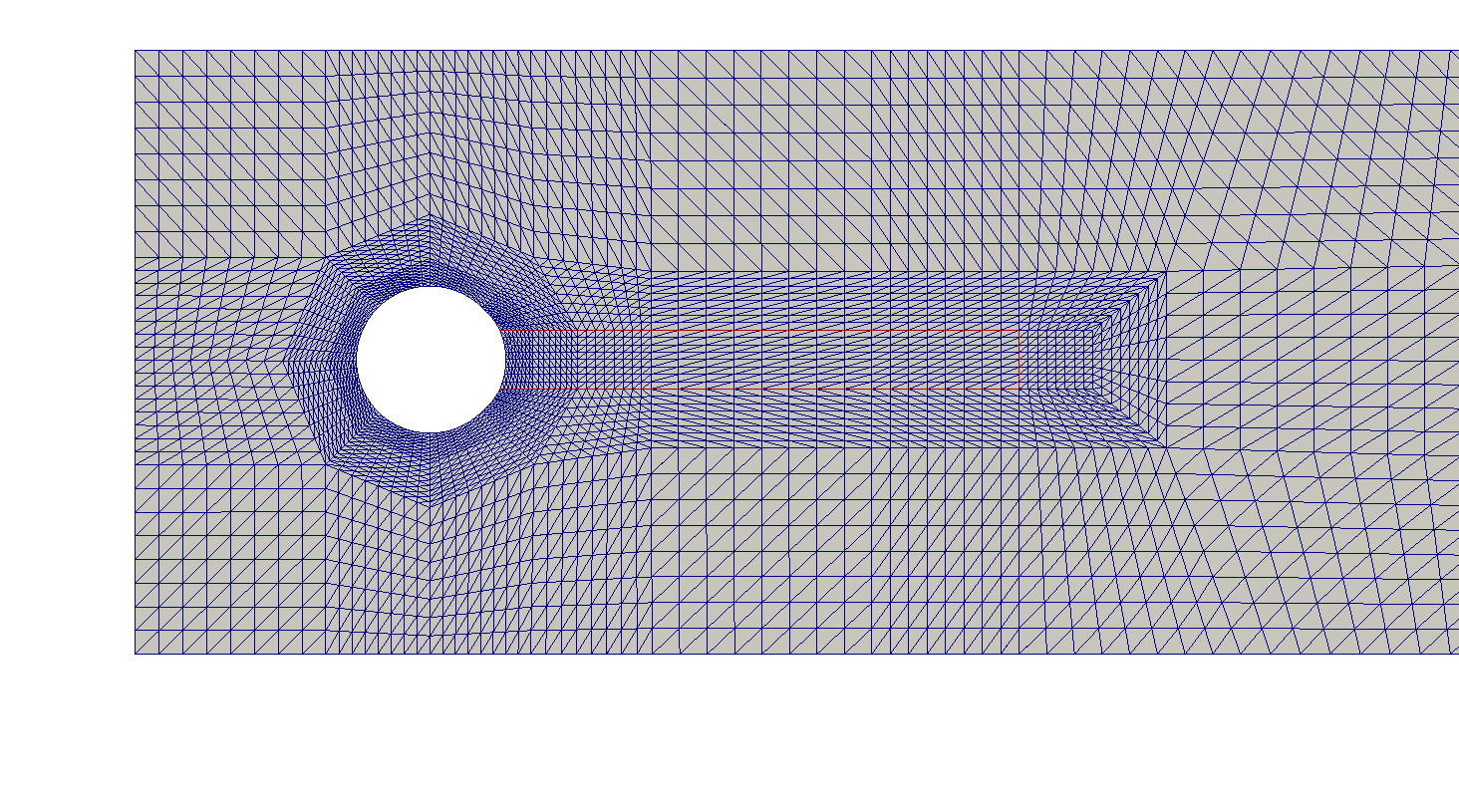}
    \caption{Example 8: Zoomed-in mesh for the fluid-structure interaction problem.}\label{fig:mesh_fsi}
\end{figure}

In the eighth example, we use the fluid-structure interaction (FSI) parameter estimation example from \cite{RiWi2013fsi}, which is an extension of the FSI-1 benchmark from Hron and Turek \cite{HronTurek2006} with a slightly thicker solid beam. The domain is shown in Figure \ref{fig:fsi_domain} and the corresponding zoomed-in mesh close to the cylindrical obstacle is shown in Figure \ref{fig:mesh_fsi}.

The domain is defined as $\Omega := (0,2.5) \times (0,0.41) \setminus B_{0.05}(0.2,0.2)$. The solid domain is $\Omega_s := \Omega \cap (0.2, 0.6) \times (0.18, 0.22)$ and the fluid domain is $\Omega_f := \Omega \setminus \Omega_s$.
The interface between the fluid and the solid domain is $\Gamma_{\text{interface}} := \Omega_s \cap \Omega_f$.  On the inflow boundary ${\Gamma_{\text{in}} = \lbrace 0 \rbrace \times (0,0.41)}$, we enforce a parabolic inflow profile, i.e. $v = g$ on $\Gamma_{\text{in}}$ with the same inflow as for the Navier-Stokes example, and for the displacement we have $u = 0$ on $\Gamma_{\text{in}}$. On the cylinder and the walls, we have homogeneous Dirichlet boundary conditions for displacement and velocity, i.e. $u = v = 0$  on $\Gamma_{\text{wall}} \cup \Gamma_{\text{circle}}$, where ${\Gamma_{\text{wall}} = \left(0, 2.5\right) \times \lbrace 0, 0.4.1 \rbrace}$ and ${\Gamma_{\text{circle}} = \partial B_{0.05}(0.2,0.2)}$.
On the boundary of the beam $\Gamma_{\text{beam}} := \partial\Omega_s \cap \partial\Omega$, we enforce homogeneous Dirichlet boundary conditions for the displacement and velocity $u = v = 0$ on $\Gamma_{\text{beam}}$. On the outflow boundary $\Gamma_{\text{out}} := \{2.5\} \times (0, 0.41)$ we use the do-nothing outflow condition $\sigma_f \cdot n = 0$ on $\Gamma_{\text{out}}$ and enforce Dirichlet boundary conditions for the displacement $u = 0$ on $\Gamma_{\text{out}}$. On the interface $\Gamma_{\text{interface}}$ we ensure that $u$ and $v$ are continuous, and we have the continuity of normal stresses $J \sigma_f F^{-T}n_f + F \Sigma_s n_s = 0$ on $\Gamma_{\text{interface}}$.
The goal of our optimization loop is to find the Lam{\'e} parameter $\mu_s$ that best matches the y-displacement at the beam tip $(0.6, 0.2)$. Once more, we define a tracking-type loss function
\begin{align*}
    J(\mu_s^{\guess}) := |u_y^{\true}(0.6, 0.2) - u_y^{\guess}(\mu_s^{\guess})(0.6, 0.2)|^2.
\end{align*}
We start with an initial guess of $\mu_s^{\guess} = 5 \cdot 10^{3}$ and try to recover the true solution $\mu_s^{\true} = 5 \cdot 10^{5}$. The other material parameters are given as $\rho_f = \rho_s = 10^3, \nu_f = 10^{-3}, \mu_f = \rho_f \nu_f, \nu_s = 0.4, \lambda_s = 2\nu_s\mu_s/ (1-2\nu_s)$. Moreover, we have
\begin{align*}
    F &= I + \nabla u, \\
    J &= \det(F), \\
    \sigma_f &= -p I + \mu_f (\nabla v F^{-1} + F^{-T}\nabla v^{T}), \\
    E_s &= \frac{1}{2}(F^T F - I) = \frac{1}{2}(\nabla u + \nabla u^T + \nabla u^T \nabla u), \\
    \Sigma_s &= 2\mu_s E_s + \lambda_s \operatorname{tr}(E_s) I.
\end{align*}
We define the function spaces
\begin{align*}
    V^u &:= \{u \in H^1(\Omega)^2 \mid u = 0 \text{ on } \Gamma_{\text{in}} \cup \Gamma_{\text{wall}} \cup \Gamma_{\text{circle}} \cup \Gamma_{\text{beam}} \cup \Gamma_{\text{out}} \}, \\
    V^v &:= \{v \in H^1(\Omega)^2 \mid u = 0 \text{ on } \Gamma_{\text{in}} \cup \Gamma_{\text{wall}} \cup \Gamma_{\text{circle}} \cup \Gamma_{\text{beam}}\}, \\
    V^p &:= L^2(\Omega) .
\end{align*}
Then, the weak formulation is given by: Find $U = (u, v, p) \in V^u \times \{(g, 0)^T + V^v\} \times V^p$ such that
\begin{align*}
    &(\rho J(\nabla v F^{-1} v), \psi^v)_{\Omega_f} + (J \sigma_f F^{-T}, \nabla \psi^v)_{\Omega_f} + (\operatorname{div}(JF^{-1}v), \psi^p)_{\Omega_f} + \alpha_u (\nabla u, \nabla \psi^u)_{\Omega_f}  \\
    &+ (F \Sigma_s, \nabla \psi^v)_{\Omega_s} + \alpha_v (v, \psi^u)_{\Omega_s} + \alpha_p [(\nabla p, \nabla \psi^p)_{\Omega_s} + (p, \psi^p)_{\Omega_s}] = 0 \\
    &\qquad \forall \Psi = (\psi^u, \psi^v, \psi^p) \in V^u \times V^v \times V^p.
\end{align*}
For the extension of the displacement to the fluid domain and of the velocity and the pressure to the solid domain, we use the parameters $\alpha_u = \alpha_p = 10^{-12}$ and $\alpha_v = 10^3$.
We use quadratic finite elements for displacement and linear finite elements for pressure to satisfy the inf-sup condition.
In our mesh, this leads to a total of 68896 degrees of freedom with 32384 DoFs each for displacement and velocity, and 4128 DoFs for pressure.

\begin{figure}[H]
    \centering
    \subfloat[Initial velocity solution]{
        \includegraphics[width=12cm]{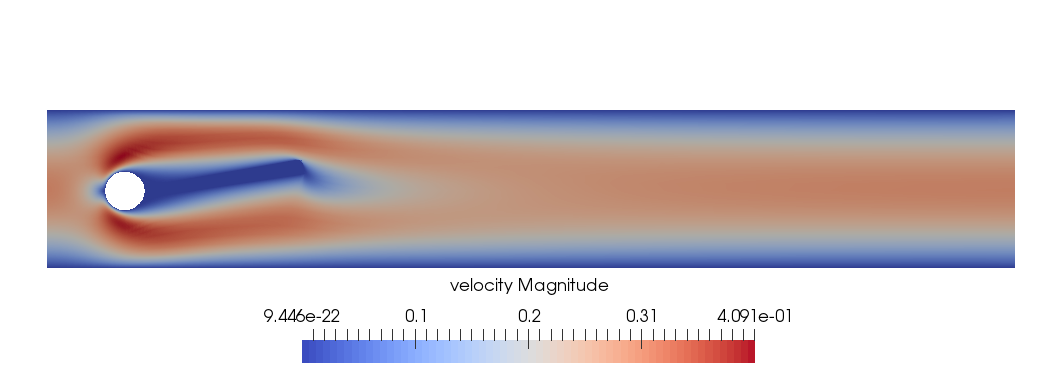}
    }\\
    \subfloat[True velocity solution]{
        \includegraphics[width=12cm]{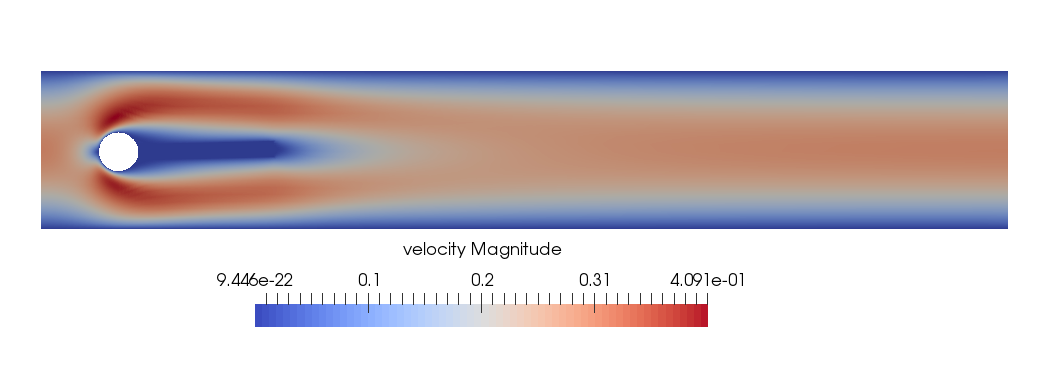}
    }
    \caption{Example 8: Initial and target solution of fluid-structure interaction with deformation scaled by a factor of 30.}\label{fig:fsi_solution}
\end{figure}

In Figure \ref{fig:fsi_solution} we plot the velocity for the initial guess ($\mu_s^{\guess} = 5 \cdot 10^3$) and for the true Lam{\'e}-parameter ($\mu_s^{\guess} = 5 \cdot 10^5$). The mesh movement has been scaled by a factor of 40 to highlight the difference between both solutions. We observe that for the initial guess we have a much larger deformation and throughout the optimization should make the material stiffer by increasing the value of $\mu_s^{\guess}$. This is also reflected in the deformation at the tip of the bar, which is given by
\begin{align*}
    u_y^{\guess}(\mu_s = 5 \cdot 10^3)(0.6, 0.2) &= 1.58370 \cdot 10^{-3}, \\
    u_y^{\true}(0.6, 0.2) = u_y^{\guess}(\mu_s = 5 \cdot 10^5)(0.6, 0.2) &= 1.80149 \cdot 10^{-4}.
\end{align*}

\begin{figure}[H]
    \centering
    \subfloat[Loss history]{
        \includegraphics[width=8cm]{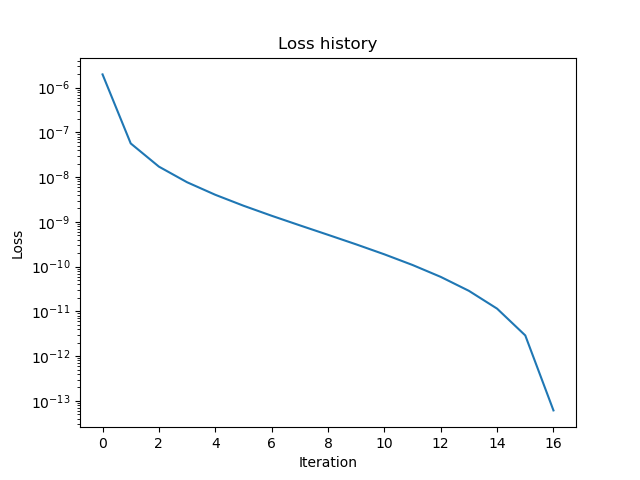}
    }
    \subfloat[Parameter history]{
        \includegraphics[width=8cm]{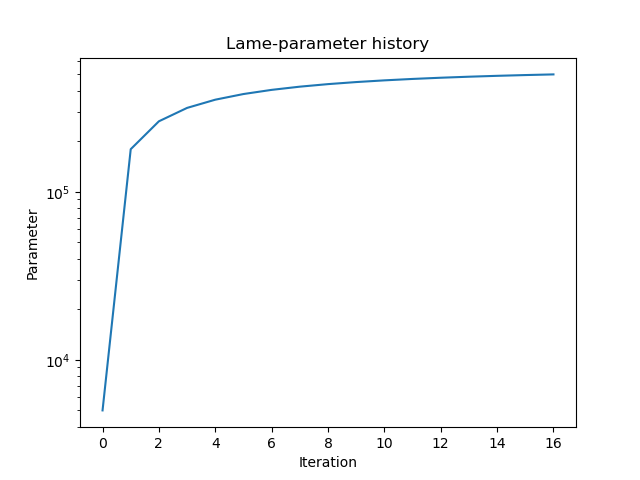}
    }
    \caption{Example 8: Loss and parameter history for fluid-structure interaction parameter estimation.}\label{fig:fsi_optimization}
\end{figure}

\noindent In Figure \ref{fig:fsi_optimization} we show the results of our parameter estimation loop where we used the Adam optimizer with a learning rate of $1.5 \cdot 10^7$.
As usual, we notice that this choice is problem-specific.
We terminate the optimization loop when an absolute tolerance of $10^{-13}$ is reached, which was the case after 16 iterations. The final guess for the Lam{\'e} parameter is $\mu_s^{\guess} =4.99225 \cdot 10^5$ which is close to the true parameter $\mu_s^{\true} =5 \cdot 10^5$. Furthermore, the final guess for the deformation of the beam is $u_y^{\guess} = 1.80397 \cdot 10^{-4}$, which is close to the true y-deformation $u_y^{\true}(0.6, 0.2) = 1.80149 \cdot 10^{-4}$, and leads to a loss of $J(4.99225 \cdot 10^5) = (u_y^{\true}(0.6, 0.2) - u_y^{\guess}(4.99225 \cdot 10^5)(0.6, 0.2))^2 = 6.13539 \cdot 10^{-14}$.

\subsection{2D Poisson with spatially-variable diffusion coefficient combined with neural networks}
\label{sec:nn_optimal_control}

In our final example, we consider a two-dimensional Poisson problem with a spatially-variable coefficient.
The strong form reads: Find $u: \Omega \rightarrow \mathbb{R}$ such that
\begin{align*}
    - \nabla \cdot (\kappa(x,y) \nabla u(x,y)) &= f, \qquad \forall (x,y) \in \Omega, \\
    u &= 0 \qquad \forall (x,y) \in \partial \Omega.
\end{align*}
Here, the right hand side is given by
\begin{align*}
    f(x, y) = -6\pi y \sin(\pi x) \cos(\pi y) + 2 \pi^2 (2x + 3y^2 + 1) \sin(\pi x) \sin(\pi y) - 2 \pi \sin(\pi y)\cos(\pi x)
\end{align*}
and the analytical solution is $u(x, y) = \sin(\pi x) \sin(\pi y)$. The true diffusion coefficient is then given as $\kappa^{\true}(x,y) = 1 + 2x + 3y^2$.
We use the loss function
\begin{align*}
    J(\kappa^{\guess}) := \|u_h^{\true} - u_h^{\guess}(\kappa^{\guess}) \|_2.
\end{align*}
In this example, we proceed differently than before and now want to find a neural network surrogate for $\kappa^{\guess}$.
This has the benefit that the optimization problem does not lead to a number of tunable parameters that is depending on the mesh elements but instead depends only on the complexity of the function $\kappa^{\true}$ and the expressivity of the neural network $\kappa_{\NN}^{\guess}: \mathbb{R}^2 \rightarrow \mathbb{R}$. For our experiments, we choose a fully connected neural network with 20 neurons in the single hidden layer and sigmoid activation function, which leads to 81 trainable parameters. We compare three approaches for this optimization problem: data-driven, physics-informed and mixed. 

In the data-driven approach we perform neural network-based regression of $\kappa^{\true}$ on a mesh with $40$ elements in x- and y-direction. In the physics-informed approach we perform only optimal control on the aforementioned mesh.
In the mixed approach we first perform regression on a coarse mesh with $4$ elements in x- and y-direction as a pre-training step and then fine-tune the neural network by solving the optimal control problem on the mesh with $40$ elements in x- and y-direction.

\begin{figure}[H]
    \centering
    \subfloat[True solution]{
        \includegraphics[width=8cm]{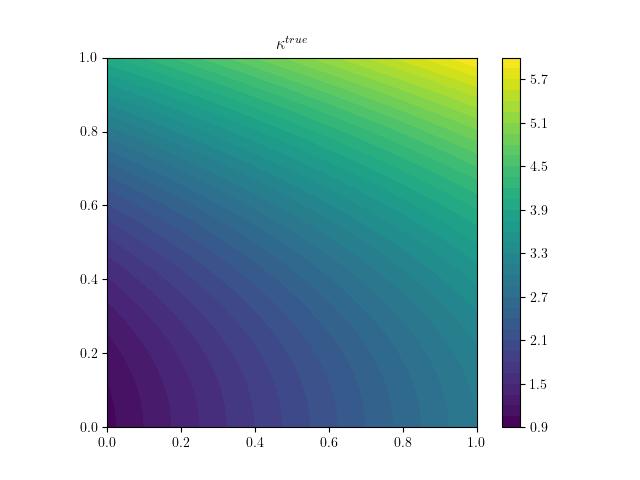}
    }%
    \subfloat[Data-driven solution]{
        \includegraphics[width=8cm]{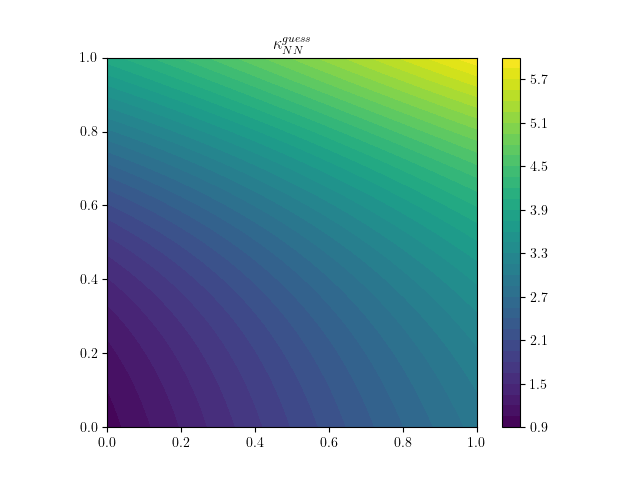}
    }\\
    \subfloat[Physics-informed solution]{
        \includegraphics[width=8cm]{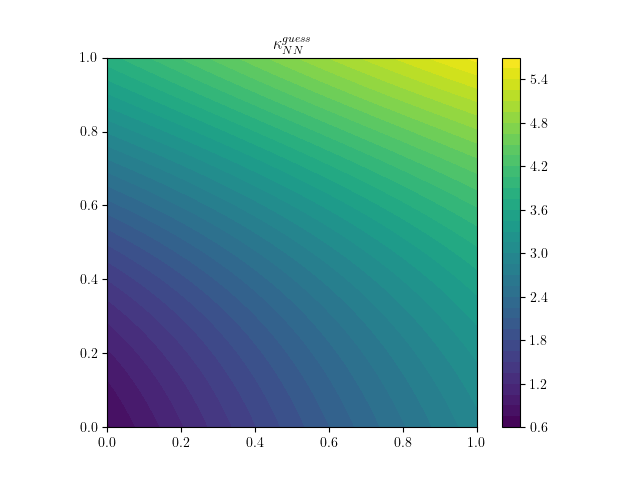}
    }%
    \subfloat[Mixed solution]{
        \includegraphics[width=8cm]{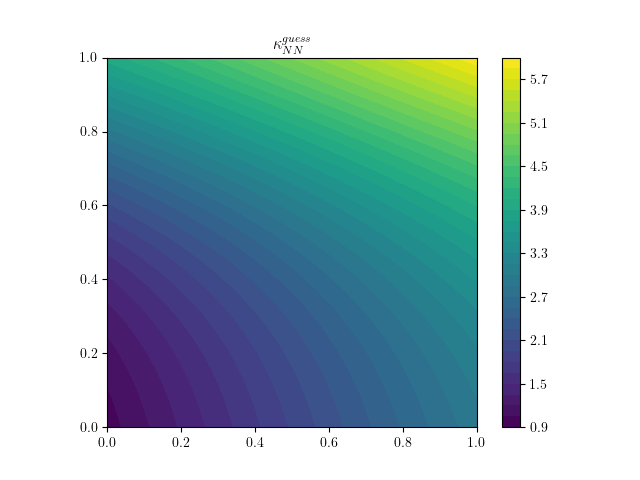}
    }%
    \caption{Example 9: Diffusion coefficient solution after neural network training.}\label{fig:poisson_NN_solution}
\end{figure}

In Figure \ref{fig:poisson_NN_solution} we see that all three approaches produce diffusion coefficients that roughly resemble the true solution.

\begin{figure}[H]
    \centering
    \subfloat[Data-driven error]{
        \hspace{-1cm}
        \includegraphics[width=6.5cm]{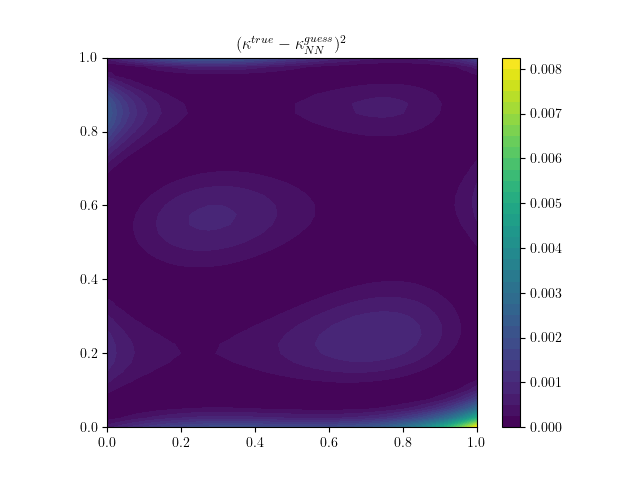}
        \hspace{-1.2cm}
    }
    \subfloat[Physics-informed error]{
        \includegraphics[width=6.5cm]{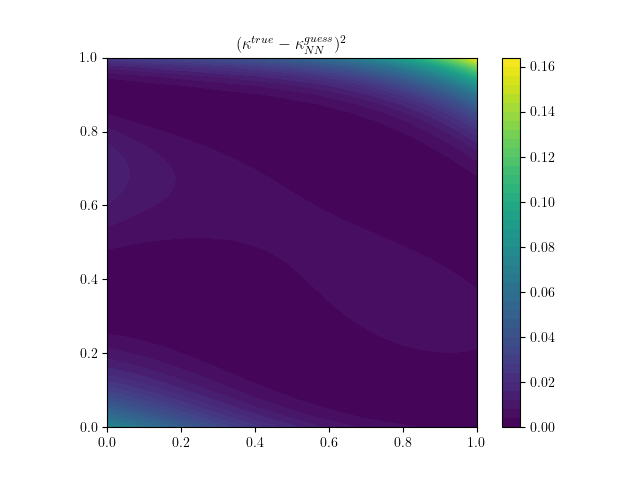}
        \hspace{-1.2cm}
    }
    \subfloat[Mixed error]{
        \includegraphics[width=6.5cm]{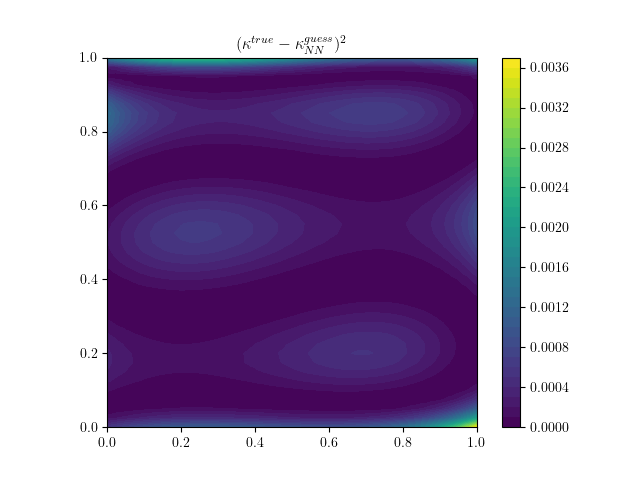}
    }%
    \caption{Example 9: Diffusion coefficient error after neural network training.}\label{fig:poisson_NN_error}
\end{figure}

In Figure \ref{fig:poisson_NN_error} we compare the squared errors between our neural network based solutions and the true diffusion coefficient $\kappa$. We observe that the purely physics-informed approach performs the worst and the data-driven error is much smaller, but also requires the knowledge of $\kappa^{\true}$ on the fine mesh. The error for the mixed approach is the lowest, but here we performed twice as many optimization steps having pre-trained on the coarse mesh and then fine-tuned on the fine mesh.

\begin{figure}[H]
    \centering
    \subfloat[Pre-trained solution]{
        \includegraphics[width=8cm]{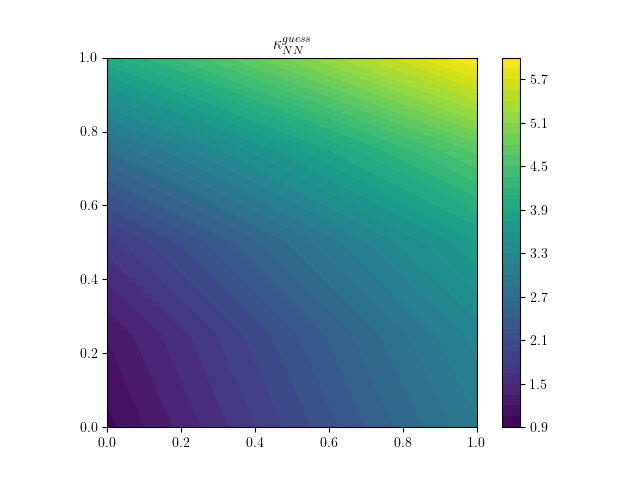}
    }%
    \subfloat[Fine-tuned solution]{
        \includegraphics[width=8cm]{finetune_solution_PoissonNN.png}
    }%
    \caption{Example 9: Diffusion coefficient error after neural network training.}\label{fig:poisson_NN_mixed}
\end{figure}

In Figure \ref{fig:poisson_NN_mixed} we show the neural network solution after pre-training on the coarse mesh in subplot (a) and after fine-tuning on the fine mesh in subplot (b). The pre-training step produces only a very rough approximation of the solution. Nonetheless, this is sufficient for a good initial guess for the physics-informed approach which then makes the training more stable and leads to a better overall model. The corresponding loss histories are presented in the figures below
\begin{figure}[H]
    \centering
        \includegraphics[width=14cm]{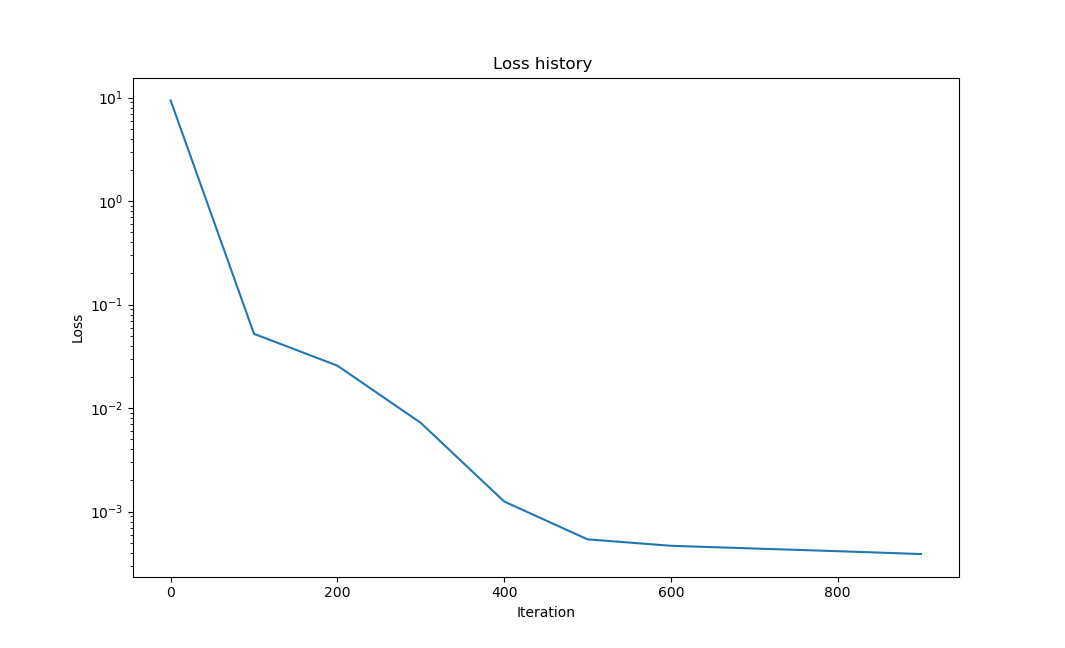}
    \caption{Example 9: Loss history for the data-driven approach.}
    \label{fig:loss_history_nn_1}
\end{figure}

\begin{figure}[H]
\centering
    \subfloat[Physics-informed loss]{
        \includegraphics[width=8cm]{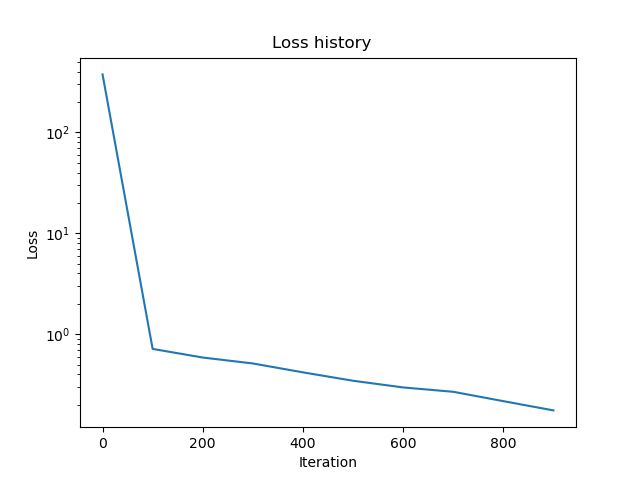}
    }%
    \subfloat[Mixed loss]{
        \includegraphics[width=8cm]{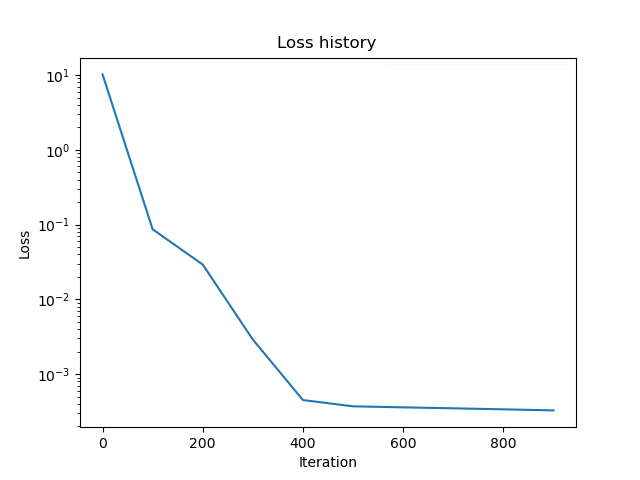}
    }%
    \caption{Example 9: Loss history for  physics-informed and mixed approach.}
    \label{fig:loss_history_nn_2}
\end{figure}

\section{Conclusion}
\label{sec:conclusion}
In this work, we gave a hands-on introduction to optimal control of PDEs in PyTorch.
To make this more accessible to a wider audience, we leveraged automatic differentiation to simplify the optimal control implementation.
By this, it is sufficient to provide a differentiable PDE solver that can be accessed from PyTorch. In the first four numerical tests, we considered linear PDEs only, i.e.,  the Poisson problem and the heat equation, using finite difference discretizations. Due to the linearity of these problems, solving the PDE is equivalent to solving a linear equation system, which was done by the differentiable direct solver (\texttt{torch\_sparse\_solve}).
For the following five numerical tests, we used finite element discretizations and considered mostly nonlinear problems, including a fluid-structure interaction problem.
Due to the nonlinearity, we did not have the equivalence between the PDE and a linear equation system anymore.
Instead, we used the automatic differentiation capabilities of \texttt{FEniCS} that have been implemented in \texttt{dolfin-adjoint}.
Therein, the backward pass for reverse mode automatic differentiation has been implemented by means of solving an auxiliary adjoint problem. 
This enabled us to solve parameter estimation, boundary control and initial condition control problems by only implementing the solver for the original PDE in \texttt{FEniCS}.
Finally, utilizing PyTorch wrappers (\texttt{torch-fenics}) allowed us to use neural network surrogates in the optimal control setting which paves the way for an easy way to combine deep learning and optimal control.

We hope that this work serves as an educational resource for people that want to learn about optimal control of PDEs and allows readers to implement their optimal control problem in only a few lines of code, especially when they want to use neural networks in this setting.

\section*{Declariations}
\subsection*{Availability of data and material}
The code for this article is available at \url{https://github.com/Denis-Khimin/OptimalControlPDEAutoDiff} 
\cite{git_repoKhiRoHeWi}.
\subsection*{Competing interests}
The authors declare that they have no known competing financial interests or personal relationships that could have appeared to influence the work reported in this paper.

\subsection*{Funding and acknowledgements}

Denis Khimin acknowledges the Priority Program 1962 (DFG SPP 1962) within the subproject
\emph{Optimizing Fracture Propagation using a Phase-Field Approach}
with the project number 314067056.
Julian Roth and Thomas Wick acknowledge the funding of the German Research Foundation (DFG)
within the framework of the International Research Training Group on
Computational Mechanics Techniques in High Dimensions GRK 2657 under Grant
Number 433082294. Alexander Henkes acknowledges support
by an ETH Zurich Postdoctoral Fellowship.


\bibliographystyle{abbrv}

\end{document}